\documentclass[10pt]{article}

\usepackage{booktabs}
\usepackage{caption}
\usepackage{times}
\usepackage{amsmath}
\allowdisplaybreaks
\usepackage{mathtools}
\usepackage{amssymb}
\usepackage{latexsym}
\usepackage{enumitem}
\usepackage{mathrsfs}
\usepackage{comment}
\usepackage{pifont}
\usepackage{bbold}
\usepackage{graphicx}
\usepackage{color}
\usepackage{hyperref}
\hypersetup{
    colorlinks=true,
    linkcolor=blue,
    filecolor=magenta,      
    urlcolor=cyan,
    citecolor=magenta
}

\newtheorem{assumption}{Assumption}

\def\qed{ \ \vrule width.2cm height.2cm depth0cm\smallskip}
\newenvironment{proof}{\noindent {\bf Proof.\/}}{$\qed$\vskip 0.1in}

\newcommand{\la}{\langle}
\newcommand{\ra}{\rangle}

\newcommand{\ol}{\overline}
\newcommand{\ul}{\underline}
\newcommand{\eps}{\varepsilon}

\newcommand{\ba}{\begin{array}}
\newcommand{\ea}{\end{array}}
\newcommand{\be}{\begin{equation}}
\newcommand{\ee}{\end{equation}}
\newcommand{\bea}{\begin{eqnarray}}
\newcommand{\eea}{\end{eqnarray}}
\newcommand{\beaa}{\begin{eqnarray*}}
\newcommand{\eeaa}{\end{eqnarray*}}

\def\dbA{\mathbb{A}}

\def\dbE{\mathbb{E}}
\def\dbF{\mathbb{F}}

\def\dbI{\mathbb{I}}

\def\dbL{\mathbb{L}}

\def\dbP{\mathbb{P}}
\def\dbR{\mathbb{R}}

\def\dbQ{\mathbb{Q}}

\def\a{\alpha}

\def\d{\delta}
\def\e{\varepsilon}

\def\l{\lambda}
\def\m{\mu}

\def\si{\sigma}
\def\t{\tau}
\def\f{\varphi}
\def\th{\theta}

\def\G{\Gamma}
\def\D{\Delta}

\def\O{\Omega}
\def\cA{{\cal A}}
\def\cB{{\cal B}}

\def\cF{{\cal F}}

\def\cJ{{\cal J}}
\def\cK{{\cal K}}
\def\cL{{\cal L}}
\def\cM{{\cal M}}
\def\cN{{\cal N}}

\def\cP{{\cal P}}

\def\cS{{\cal S}}

\def\cW{{\cal W}}

\def\no{\noindent}

\def\ms{\medskip}

\def\q{\quad}
\def\qq{\qquad}

\def\pa{\partial}
\def\cd{\cdot}
\def\cds{\cdots}

\def\tr{\hbox{\rm tr}}

\def\qed{ \hfill \vrule width.25cm height.25cm depth0cm\smallskip}

\newcommand{\basa}{\begin{assumption}}
\newcommand{\easa}{\end{assumption}}

\newcommand{\bas}{\begin{assum}}
\newcommand{\eas}{\end{assum}}

\def\pa{\partial}

 \def\cd{\cdot}
\def\cds{\cdots}

\def\tr{\hbox{\rm tr$\,$}}

\def\dis{\displaystyle}

\def\bQ{{\bf Q}}

\def\1{{\bf 1}}

\def\:{\!:\!}
\def\reff#1{{\rm(\ref{#1})}}
\def \proof{{\noindent \bf Proof\quad}}

\def \bS{{\bf S}}
\def \bm{{\bf m}}

at 9pt

\begin{document}

\newtheorem{thm}{Theorem}[section]
\newtheorem{lem}[thm]{Lemma}
\newtheorem{cor}[thm]{Corollary}
\newtheorem{prop}[thm]{Proposition}
\newtheorem{rem}[thm]{Remark}
\newtheorem{eg}[thm]{Example}
\newtheorem{defn}[thm]{Definition}
\newtheorem{assum}[thm]{Assumption}

\numberwithin{equation}{section}

\title{\bf{Viscosity solutions for mean field optimal switching with a two-time-scale Markov chain\footnote{The authors are grateful for the financial support from the Natural Science Foundation of China (11831010, 61961160732), the National Key Research and Development Program of China (2023YFA1009200), Shandong Province Natural Science Foundation, China (ZR2019ZD42) and the Taishan Scholars Climbing Program of Shandong (TSPD20210302).}}}
\author{Tian Chen\footnote{Shandong University, Jinan, China,  chentian43@sdu.edu.cn} \quad Guanxu Li\footnote{Shandong University, Jinan, China, lgx96@mail.sdu.edu.cn} \quad Zhen Wu\footnote{Shandong University, Jinan, China, wuzhen@sdu.edu.cn }}
\date{}

\maketitle

\begin{abstract}
In this paper, we consider the mean field optimal switching problem with a Markov chain under viscosity solution notion. Based on the conditional distribution of the Markov chain, the value function and corresponding dynamic programming principle (DPP) are established. The switching problem is characterized by an obstacle equation on the Wasserstein space, and the existence, stability, and comparison principle are obtained in the sense of viscosity solution. In particular, we consider a two-time-scale structure and obtain the convergence of the limit system. As an application of our theoretical results, an innovative example concerning the stock trading problem in a regime switching market is solved.
\end{abstract}

\no{\bf MSC2020.} 60J10, 35Q89,  49N80, 49L25, 49J40.

\vspace{3mm}
\no{\bf Keywords.} Mean field optimal switching,  viscosity solutions, obstacle problem, Markov chain, two-time-scale.

\section{Introduction}

The viscosity solution theory is a celebrated and powerful tool in optimal control problems. Specifically, the wellposedness of the dynamic programming equation can be proved under weaker regularity for the smoothness of the value function is not satisfied in general. The pioneering works of this topic may refer to \cite{CEL,CIL,CL}. In recent years, viscosity solution theory has been further developed with the rise of mean field problems. For example, the McKean-Vlasov control can be characterized by master equations, the wellposedness of which relies on viscosity solution theory on the Wasserstein space, \cite{WZ}. Along this direction, this paper is focus on the mean field optimal switching problem with a Markov chain under viscosity solution notion. By virtue of the conditional distribution of the Markov chain,  we characterize the problem by an obstacle problem on the Wasserstein space through dynamic programming principle (DPP).

There have been some serious efforts on viscosity solutions of the nonlinear partial differential equations on Wasserstein space. One natural idea is to modify the definition of sup/sub differential to adapt the setting of Wasserstein space. Along this line, Cardaliaguet \& Quincampoix \cite{CarQui} considered a first order Hamilton-Jacobi-Isaacs equation on Wasserstein space, which is arose from deterministic zero-sum games with random initial conditions. The comparison principle for viscosity solutions was established by combining the doubling variables argument with Ekeland's variational principle. We may also mention the work of Gangbo, Nguyen \& Tudorascu \cite{GNT} and Jimenez, Marigonda \& Quincampoix \cite{JMQ}, who also define a notion of viscosity solutions for Hamilton-Jacobi equations by using sub-differentials. Another approach followed by several authors (see e.g. \cite{PW}) consists in exploiting Lions' idea \cite{Lions} of lifting the functions on Wasserstein space into functions on Hilbert space of random variables and then using the existing viscosity theory on Hilbert spaces (see e.g. \cite{FGS,Lions1,Lions2,Lions3}). More recently, Cosso et al. \cite{CGKPR} defined viscosity solutions for Hamilton-Jacobi equations by requiring the global extrema on the Wasserstein space for the tangency property of the test functions.

In the context of mean field control problems in a path dependent setting, Wu \& Zhang \cite{WZ} proposed a notion of viscosity solutions for parabolic equations on the Wasserstein space by restricting the viscosity neighborhood of some point $(t,\mu)$ (where $t$ is a time and $\mu$ a measure) on the Wasserstein space to ensure the compactness. More recently, Talbi, Touzi \& Zhang \cite{TTZ2} studied the mean-field stopping problem employing the same notion of viscosity solution. We would also like to mention Soner \& Yan \cite{SY}, in which the Sobolev norms on the space of measures are  adopted to guarantee the compactness of the viscosity neighborhood.
We shall follow the approach of \cite{WZ} and  consider the joint law of $(X_t, I_t)$ as the variable of the value function, where $X$ is the state process and $I$ is the switching process. The presence of Markov chain $\a_t$ makes things more subtle so that we need consider the conditional distribution  to construct an appropriate value function. Some new development related to DPP contained Markov chain can refer to \cite{GW,LV}. We show that, under natural conditions, the value function of the mean field optimal switching problem is indeed the unique viscosity solution of the corresponding obstacle equation on Wasserstein space.  It should be pointed that the appearance of Markov chain brings more computational complexity. To overcome the difficulty, we employ the so-called two-time-scale approach  and obtain a convergence result for associated limit problem in Section \ref{sect-TTS}.

The main contributions of this paper include the following:

(i) We formulate a general mean field optimal switching problem under a regime-switching model, and establish the corresponding DPP by considering the conditional distribution of the Markov chain, which can be seen as a generalization of the standard case. (ii) We prove that the value function of the mean field optimal switching problem is indeed the unique viscosity solution of the corresponding obstacle equation on Wasserstein space. (iii) When the Markov chain has a two-time-scale structure, the viscosity solution approach is employed to verify the convergence of the value function, and thus help us to simplify the computational procedure. (iv) The mean field optimal switching framework specially provides an appropriate and general model to study the massive stock trading problem. The optimal strategy is obtained in a regime switching market, which is characterized by the Markov chain.

The remaining sections are organized as follows. In Section \ref{sect-obstacle}, we present the mean field optimal switching problem with regime-switching, the corresponding DPP and variational inequality. We propose our definition of viscosity solutions and prove the main results in Section \ref{sect-viscosity}. Section \ref{sect-TTS} is devoted to two-time-scale approach, which is used for the stock trading problem in a regime switching market. Section \ref{sect-conc} concludes the paper. Finally, we prove some technical results in the Appendix.

\textbf{Notations}. We denote by $\cP(\O,\cF)$ the set of probability measures on a measurable  space $(\O,\cF)$, and $\cP_2(\O,\cF) := \{ m \in \cP(\O,\cF) : \int_\O d(x_0,x)^2 m(dx)$ $ < \infty \}$ for some $x_0 \in \O$, where $d$ is a metric on $\O$. $\cP_2(\O,\cF)$ is equipped with the corresponding $2$-Wasserstein distance $\cW_2$. When $(\O,\cF) = (\dbR^d,\cB(\dbR^d))$, we simply denote them as $\cP(\dbR^d)$ and $\cP_2(\dbR^d)$. For a random variable $Z$ and a probability $\dbP$, we denote by $\dbP_Z:=\dbP\circ Z^{-1}$ the law of $Z$ under $\dbP$.
For vectors $x, y\in \dbR^n$ and matrices $A, B\in \dbR^{n\times m}$, denote $ x\cd y:=\sum_{i=1}^n x_iy_i$  and $A:B:= \tr(A B^\top)$. We shall also write USC (resp. LSC) upper (resp. lower) semi-continuous".

\section{The obstacle problem on Wasserstein space}\label{sect-obstacle}

\subsection{Formulation}
Let $T < \infty$ be fixed, and $\O := C^0([-1,T],\dbR^d) \times \dbI^0([-1,T])\times \dbA^0([-1,T])$ the canonical space, where: 

\no$\bullet$ $C^0([-1,T],\dbR^d)$ is the set of continuous paths from $[-1,T]$ to $\dbR^d$, constant on $[-1,0)$; \\
\no$\bullet$ $\dbI^0([-1,T])$ is the set of c\`adl\`ag maps from $[-1,T]$ to $\cN$,  constant  on $[-1,0)$;\\
\no$\bullet$ $\dbA^0([-1,T])$ is the set of c\`adl\`ag maps from $[-1,T]$ to $\cM$,  constant  on $[-1,0)$.
 
 \no We equip $\O$ with the Skorokhod distance, under which it is a Polish space. The choice of the extension to $-1$ is arbitrary,  the extension is only needed to allow for an immediate switching at time $t=0$.

Define switching process $I_t :=\underset{n \geq 1}{\sum}\xi_{n-1}\1_t^n$, where $\1_t^n=\1_{\t_{n-1} \leq t < \t_n}$ with  $\t_0:=-1$. Let  $\xi_n $ take value in a finite set $\cN$, and Markov chain $\a_t$ take value in a finite set $\cM$. Denote $Q = [\l_{pq}(x)]$ the generator of $\a$ with $\lambda_{pq}\ge 0$ for $p \neq q$ and $\underset{q \in \cM}{\sum}\lambda_{pq}=0$,  $\forall p \in \cM$. The evolution of the controlled system is governed by:
\bea
\label{dynamic}
d X_t=b(\a_t,X_t,I_t,\cL_{X_t})dt+\si(\a_t,X_t,I_t,\cL_{X_t})dW_t,\q X_0=x,
\eea  
here $\cL_{X_{\cdot}}$ is the law of the state process, which serves as a mean field term in our structure. We denote $Y:=(X,I)$ and $\a$ the canonical process, with corresponding state space $\mathbf{S}:=\dbR^d\times\cN$ and $\cM$, canonical filtration $\dbF = (\cF_t)_{t \in [-1,T]}$, and the corresponding switching time of the process $I$:
\bea
\label{tau}
\t_n := \inf\{t > \t_{n-1} : I_t \neq I_{t-}\},  \ \mbox{for all $t \in [0,T)$.}
\eea 
By the c\`adl\`ag property of $I$, $\t_n$ is an $\dbF-$stopping time. Our problem consists in maximizing over the switching controls $(\t_n,\xi_n)$ the gain functional 
\bea
\label{functional}
J_0:=\dbE\Big[ h(\a_T,X_T,I_T,\cL_{X_T})+\int_0^T f(s,\a_s, X_s,I_s, \cL_{X_s})ds-\sum_{n \geq 1}g_{\xi_{n-1}\xi_n}(\t_n,X_{\t_n})\mid \a_0=q \Big] .
\eea 
Denote further 
 \beaa
 \mathbf{Q}_t:=[t,T)\times\cP_2(\mathbf{S}\times \cM), 
 &\mbox{and}&
 \overline{\mathbf{Q}}_t:=[t,T]\times\cP_2(\mathbf{S}\times \cM),
 ~~t\in[0,T).
 \eeaa

 Let $(b,\si,f): [0,T] \times \mathbf{S} \times \cM \times \cP_2(\mathbf{S}) \rightarrow \dbR^d \times \cS_d^+ \times \dbR$ and $h:\cM \times \mathbf{S} \times \cP_2(\mathbf{S}) \rightarrow \dbR, g : [0,T] \times \dbR^d \rightarrow \dbR$, where  $\cS_d^+$ denotes the set of $d\times d$  non-negative symmetric matrices. We introduce the joint marginal law $m_t := \dbP_{(X_t,I_t)}$ and denote $m(i)=\dbP(I_t=i)$. Throughout the paper, the following assumption will always be in force,  where $\cP_2$ is equipped with the $\cW_2$-distance. 
 
\begin{assum}
\label{assum-bsig}
\no{\rm (i)} $b, \si$ are continuous in $t$, and uniformly Lipschitz continuous in $(q,x, m)$. \\
{\rm (ii)}  $f$ is Borel measurable and has quadratic growth in $x\in \dbR^d$, and the following function $F$ is continuous on $\overline{\mathbf{Q}}_0$: 
\bea
\label{F}
F(t,m,l) := \sum_{q \in \cM}\int_{\bS} f(t,q,y,m)m(dy)l(q).
\eea
{\rm (iii)} $h$ is locally bounded and has quadratic growth in $x\in \dbR^d$; and extended to $\cP_2(\bS)$ by $h(\cdot,y,m) := h(\cdot,y,m(\cdot, \cN))$. Suppose the following function $H$ is continuous on $\cP_2(\bS\times\cM)$: 
\bea
\label{H}
H(m,l) := \sum_{q \in \cM}\int_{\bS} h(q,y,m)m(dy)l(q).
\eea
{\rm (iv)} Denote $Q(m)= [\l_{pq}(m)]$ the extended generator of $\a$,  and the following is uniformly bounded:
\bea
\label{Q}
\l_{pq}(m) := \underset{i \in \cN}{\sum}\int_{\dbR^d} \l_{pq}(x)m(dx,i).
\eea
{\rm (v)}  $g$ has quadratic growth in $x\in \dbR^d$, and the following function $G^n$ is continuous on $[0, T]\times \cP_2(\bS)$:
\bea
\label{G}
G^n(t,m) := \underset{i \in \cN}{\sum}\underset{j \in \cN}{\sum}\int_{\dbR^d} g_{ij}(t,x)m^{i}(dx,j)\dbP(\xi_{n-1}=i),
\eea
where $m^{i}(\cdot,j)$ is the conditional probability of switching from the previous state $i$ to current state $j$.
 Similar to the standard case, for given $(t,x)\in [0, T]\times \dbR^d$,
\bea
\label{g}
g_{ij}(t,x)+g_{jk}(t,x)>g_{ik}(t,x),\q g_{ii}(t,x)=0,\q \forall i,j,k \in \cN.
\eea 
\end{assum}

We define the McKean-Vlasov dynamics on $[0, T]$ in the sense of weak formulation:
\bea
\label{WP}
X_s = X_0 + \int_0^s b(\a^{\dbP}_r, Y_r, \dbP_{Y_r})dr + \int_0^s \si(\a_r^{\dbP}, Y_r, \dbP_{Y_r})dW_r^\dbP ,
\eea
where a solution $\dbP$ of the last SDE is defined by the requirement that the following processes $M$ and $N$ are $\dbP-$martingales on $[0, T]$:
 \bea\label{martingalepb}
M_. := X_. - \int_0^. b(\a^{\dbP}_r,Y_r,\dbP_{Y_r})dr \ \mbox{and} \ N_.:= M_.^2 -  \int_0^. \si^2(\a^{\dbP}_r,Y_r,\dbP_{Y_r})dr.
\eea 
Here $\a^{\dbP}$ is a Markov chain under $\dbP$ with  extended generator $Q(\dbP_Y)$. It should be noted that $\a^{\dbP}$ is uniquely identified by $\dbP$, see Douc et al.\cite[Vol. I, Theorem 3.1.2]{DMPS}.
 We then focus on the mean field optimal switching problem: given $(m,q) \in \cP_2(\bS) \times \cM$,
\beaa\label{naive}
V_0 := \sup_{\dbP} \dbE^\dbP\Big[h(\a^\dbP_T,Y_T,\dbP_{Y_T})+\int_0^T  f(r,\a^\dbP_r, Y_r, \dbP_{Y_r})dr-\sum_{\t^{\dbP}_n \geq 0}g_{\xi_{n-1}\xi_n}(\t^\dbP_n,X_{\t^\dbP_n})\mid \dbP_{\a_{0-}}=\d_q \Big] \\*
= \sup_{\dbP} \Big[H(\dbP_{\a_T^{0,\d_q}},\dbP_{Y_T})+\int_0^T F(r,\dbP_{\a_r^{0,\d_q}}, \dbP_{Y_r})dr-\sum_{\t^{\dbP}_n \geq 0}\int_0^T G^n(r, \dbP_{Y_r})dr\Big]. 
\eeaa
where the supremum is taken over all solutions $\dbP$ of the McKean-Vlasov SDE satisfying the constraint $\dbP_{Y_{0-}} = m$ and $\dbP_{\a_{0-}} = \d_q$. Here $\a^{t,l}$ denotes the Markov chain starting from $(t,l)\in [0,T]\times \cP_2(\cM)$.

\subsection{Dynamic programming principle}
 In order to solve this problem, we use the dynamic programming approach by taking the  distribution $m_t$, $l_t$ of the variables $Y_t=(X_t,I_t)$, $\a_t$ respectively as  new variables, which leads to the dynamic value function by considering conditional distribution on Markov chain.
\begin{equation}\label{weakoptstop}
 V (t,m;l) := \underset{\dbP \in \cP(t,m;l)}{\sup} \Big[H(\dbP_{\a^{t,l}_T},\dbP_{Y_T})+\int_t^T F(r, \dbP_{\a^{t,l}_r}, \dbP_{Y_r})dr-\sum_{\t^{\dbP}_n \geq t}\int_t^T G^n(r, \dbP_{Y_r})dr\Big], \q  \mbox{$(t,m,l) \in  \bQ_0$},
\end{equation}
where $\cP(t,m;l)$ is the set of probability measures $\dbP$ on $(\O,\cF_T)$ such that 

\no$\bullet$ $\dbP_{Y_{t-}} = m$ and $s\in [-1, t) \to Y_s$ is constant, $\dbP$-a.s.

\no$\bullet$ $\dbP_{\a_{t-}} = l$ and extended generator equals $Q(\dbP_{Y_s})$, for $s\in [t,T]$, $\dbP$-a.s.

\no$\bullet$ The processes $M,N$ of \eqref{martingalepb} are $\dbP$-martingales on $[t,T]$, so that, for some $\dbP$-Brownian motion $W^\dbP$,
\begin{equation}\label{asympt}
X_s = X_t + \int_t^s b(\a^{\dbP}_r, Y_r, \dbP_{Y_r}) dr + \sigma(\a^{\dbP}_r, Y_r, \dbP_{Y_r})dW_r^\dbP , \ \dbP-\mbox{a.s.}
\end{equation}

 A special element of $\cP(t,m;l)$ is $\bar \dbP$ under which $X$ is unstopped. That is,
\bea
\label{barP}
X_s = X_t + \int_t^s b(\a^{\bar \dbP}_r, Y_r, \bar \dbP_{Y_r})dr + \int_t^s \si(\a^{\bar \dbP}_r, Y_r, \bar \dbP_{Y_r}) dW^{\bar \dbP}_r,~ I_s = I_{t-}\1_{[t, T)}(s),~ \bar \dbP \mbox{-a.s.} 
\eea

\begin{lem}\label{existence} 
For any $(t,m,l)\in \bQ_0$, the set $\cP(t,m;l)$ is compact under the Wasserstein distance $\cW_2$. 
\end{lem}
The above lemma guarantees the existence of optimality since $F,H,G$ is continuous. The proof is slightly different with Proposition 2.2 in Talbi, Touzi \& Zhang \cite{TTZ}, we sketch it here.

\proof
We assume for simplicity $t=0$. Firstly, we can prove the uniform integrability following the same procedure in \cite{TTZ}, i.e. 
\bea
\label{qUIlem}
\sup_{\dbP\in \cP(0,m;l)}\dbE^\dbP\big[ |X^*_T|^2\big]\le C_m,\qq \lim_{R \rightarrow \infty}\sup_{\dbP\in \cP(0,m;l)}\dbE^{\dbP}\big[ |X^*_T|^2 \1_{\{X^*_T \ge R\}}\big] = 0.
\eea
Here $X^*_T := \sup_{0\le s \le T} |X_s|$.

We next show $\cP(0,m;l)$ is closed under weak convergence. Let $\{\dbP^n\}_{n \ge 1} \subset \cP(0,m)$ converge weakly to some $\dbP^\infty$. Since $\dbP^n_{(X_0,I_{0-})} = m$  for all $n$, we have  $\dbP^\infty_{Y_{0-}} = m$. Similarly, we have $\dbP^\infty_{\a_{0-}} = l$ and the extended generator equals $Q(\dbP^\infty_{Y})$. Then it suffices to show that the processes $M,N$ in \eqref{martingalepb} are $\dbP^\infty-$martingales on $[0,T]$. We shall report only the detailed argument for $M$, as it is immediately adapted to $N$. 

By the Skorokhod's representation theorem(see \cite{B2013}), there exists a probability space $(\O^0, \cF^0, \dbP^0)$ and processes $\{(Y^n, \a^n)\}_{n \ge 1}$ and $(Y^\infty, \a^\infty)$ defined on this space such that,  
 \bea
 \label{XInconv}
 \dbP^n_{(Y,\a)} = \dbP^0_{(Y^n,\a^n)} \ \mbox{for all $n \le \infty$, and} (Y^n, \a^n) \overset{d_{SK}}{\longrightarrow} (Y^\infty,\a^\infty), \ \dbP^0-\mbox{a.s.}
 \eea
 For all $n \ge 1$, the $\dbP^n-$martingale property of $M$ translates to:
 \bea\label{projsko}
  \dbE^{\dbP^0}[(M_s^n-M_{t}^n)\psi(Y_{. \wedge {t}}^n)] = 0 \q \mbox{ for all $\psi \in C_b(\O)$ and $0\le t\le s\le T$,}
 \eea
 with $M^n_s = X^n_{t} - \int_{t}^s b(\a^{n,\dbP^0}_r,Y^n_r,\dbP^0_{Y^n_r})dr$ and $C_b(\O)$ the set of $\dbR^d$-valued bounded continuous functions on $\O$. Moreover,  for $r \in [t,T]$, by the Lipschitz continuity of $b$ we have
 \beaa
 \lvert b(\a^{n,\dbP^0}_r,X_r^n, \dbP^0_{Y_r^n}) - b(\a^{\infty,\dbP^0}_r,X_r^\infty, \dbP^0_{Y_r^\infty})\rvert &\le& C\big[d_{SK}(\a^{n,\dbP^0}_r,\a^{\infty,\dbP^0}_r)+\lvert X_r^n - X_r^\infty \rvert + \cW_2(\dbP^0_{Y_r^n},\dbP^0_{Y_r^\infty})\big].
 \eeaa
Send $n\to \infty$, by \reff{XInconv} we have $\lvert X_r^n - X_r^\infty \rvert \to 0$, $\dbP^0$-a.s. and 
  $\cW_2(\dbP^0_{Y_r^n},\dbP^0_{Y_r^\infty}) \to 0$.
  Thus
  $$  b(\a^{n,\dbP^0}_r,X_r^n, \dbP^0_{Y_r^n}) \underset{n \rightarrow \infty}{\longrightarrow} b(\a^{\infty,\dbP^0}_r,X_r^\infty, \dbP^0_{Y_r^\infty}), \ \mbox{$\dbP^0$-a.s.} $$
  Moreover, as $b$ is Lipschitz and $\{X^n\}_{n \ge 1}$ are uniformly integrable (as the 2-uniform integrability of \reff{qUIlem} implies the 1-uniform integrability), then $\{M^n\}_{n\ge 1}$ are uniformly integrable. The convergence for the Skorokhod distance also implies the convergence of $I_{. \wedge t}^n$ to $I_{. \wedge t}^\infty$.
This allows to take the limit in \eqref{projsko} as $\psi \in C_b(\O)$, hence
 $ \dbE^{\dbP^0}[(M_s^\infty-M_t^\infty)\psi(Y_{. \wedge t}^\infty)] = 0.$
By the arbitrariness of $\psi \in C_b(\O)$, this proves $M^\infty$ is a $\dbP^0-$martingale, or equivalently that $M$ is a $\dbP^\infty-$martingale. 

Finally, we show that $\cP(0, m;l)$ is compact under $\cW_2$. Let $\{\dbP^n\}_{n\ge 1}\subset \cP(0, m;l)$. First, one can easily obtain a uniform bound for the conditional variation of $Y$ under all $\dbP^n$, then by Meyer \& Zheng \cite[Theorem 4]{ZM} we see that $\{\dbP^n\}_{n\ge 1}$ is relatively weakly compact, namely there exists a weakly convergent subsequence.  By closed $\cP(0,m;l)$,  without loss of generality we assume the whole sequence $\dbP^n\to \dbP^\infty\in \cP(0, m;l)$ weakly. Then it follows from Carmona \& Delarue \cite[Vol. I, Theorem 5.5]{CarDel} that $\underset{n \rightarrow \infty}{\lim}\cW_2(\dbP^n,\dbP^\infty) = 0$. This proves the compactness of $\cP(0, m;l)$.
\qed

Hence, the DPP can be constructed similar with \cite{LW}. We only state the difference in the proof.
\begin{prop}\label{DPP}
 For any $(t,m,l)\in \bQ_0$ and $s\in [t, T]$, we have the DPP:
\bea\label{weakDPP}
 V(t,m;l) = \sup_{\dbP \in \cP(t,m;l)} \Big[\int_t^s F(r, \dbP_{\a^{t,l}_r},\dbP_{Y_r})dr -\underset{\t^{\dbP}_n \geq t}{\sum}\int_t^s G^n(r, \dbP_{Y_r})dr + V(s, {\dbP}_{Y_{s-}},\dbP_{\a^{t,l}_{s-}})\Big] .
\eea 
\end{prop}

\proof
Denote, for any probability measure $\dbP$ on $(\O, \cF_T)$,
  \begin{eqnarray*}
  J(t, \dbP;l) &:= H(\dbP_{\a^{t,l}_T},\dbP_{Y_T})+\int_t^T F(r,\dbP_{\a^{t,l}_r}, \dbP_{Y_r})dr-\underset{\t^{\dbP}_n \geq t}{\sum}\int_t^T G^n(r, \dbP_{Y_r})dr.
 \end{eqnarray*}
 Let $\tilde V(t, m;l)$ denote the right term of (\ref{weakDPP}). For arbitrary $\dbP\in \cP(t, m;l)$, denote $\tilde m:= \dbP_{Y_{s-}}$, $\tilde l:= \dbP_{\a^{t,l}_{s-}}$.
  
  First, we construct the finite measure by adding the distribution of Markov chain as follows. For any $A_i\in \cB(\bS)$ and $B_i\in \cB(\cM)$,
\small{
\begin{eqnarray*}
\nu_\pi((A_0,B_0)\times \cds \times (A_{m+n},B_{m+n})) := \dbP\Big(Y_{s-} \in \overset{m}{\underset{i=0}{\cap}} A_i, ~\a_{s-}\in B_m,~ (Y_{t_{m+j}},\a_{t_{m+j}}) \in (A_{m+j},B_{m+j}), j=1,\cds, n\Big). 
\end{eqnarray*}
}
Unlike $Y_{s-}$, $\a_{s-}$ is only required $B_m$-measurable since the Markov property. Then it follows from the Kolmogorov extension theorem that there exists a probability measure $\tilde\dbP$ on $(\O, \cF_T)$ such that $\{\nu_\pi\}_\pi$ is the finite distribution of the process $(Y,\a)$ under $\tilde\dbP$. It is straightforward to verify $\tilde\dbP\in \cP(s, \tilde m;\tilde l)$, and $\tilde\dbP_{(Y_r,\a_r)} = \dbP_{(Y_r,\a_r)}$ for all $r\in [s, T]$. Then it's clear that 
$$J(t, \dbP;l)\leq \int_t^s F(r,\dbP_{\a^{t,l}_r}, \dbP_{Y_r})dr-\sum_{\t^{\dbP}_n \geq t}\int_t^s G^n(r, \dbP_{Y_r})dr + V(s, {\dbP}_{Y_{s-}};\dbP_{\a^{t,l}_{s-}}).$$
 Since $\dbP\in \cP(t,m;l)$ is arbitrary, we obtain $ V(t,m;l) \le \tilde V(t,m;l) $. 

On the other hand,  $\forall \tilde m$ and $\tilde l$, by Lemma \ref{existence} there exists $\tilde\dbP\in \cP(s, \tilde m;\tilde l)$ such that $J(s, \tilde \dbP;\tilde l) = V(s, \tilde m;\tilde l)$. For the above time partition $\pi$, we introduce another finite measure: for any $A_i\in \cB(\bS)$,
\small{
\begin{eqnarray*}
\nu_\pi((A_0,B_0)\times \cds \times (A_{m+n},B_{m+n})) &:=& \sum_{q\in \cM}\tilde l(q)\int_{\bS}\dbE^\dbP\Big[  \prod_{i=0}^m \1_{A_i}(Y_{t_i}) \prod_{i=0}^m \1_{B_i}(\a_{t_i}) \Big|Y_{s-} = y,\a_{s-} = q\Big]\\
&\times& \dbE^{\tilde \dbP}\Big[ \prod_{j=1}^n \1_{A_{m+j}}(Y_{t_{m+j}})\prod_{j=1}^n \1_{B_{m+j}}(\a_{t_{m+j}}) \big| Y_{s-}=y,\a_{s-} = q\Big] \tilde m(dy). 
\end{eqnarray*}
}
Applying the Kolmogorov extension theorem again there exists a probability measure $\hat\dbP$ on $(\O, \cF_T)$ such that $\{\nu_\pi\}_\pi$ is the finite distribution of the process $(Y.\a)$ under $\hat\dbP$. It is clear that $\hat\dbP = \dbP$ on $\cF_{s-}$, and $\{(Y_{s-},\a_{s-}), (Y_r,\a_r), s\le r\le T\}$ has the same distribution under $\hat \dbP$ and $\tilde \dbP$. Note that the presence of $\a$ has no influence on the conditional independence under $\hat \dbP$ on $(Y_{s-},\a_{s-})$, we can check that $M,N$ are $\hat \dbP$-martingales on $[t, T]$ same with \cite{LW}. Therefore, by setting $(\tilde m. \tilde l)=({\dbP}_{Y_{s-}};\dbP_{\a^{t,l}_{s-}})$
\begin{eqnarray*}
&&\dis \int_t^s F(r,\dbP_{\a^{t,l}_r}, \dbP_{Y_r})dr-\sum_{\t^{\dbP}_n \geq t}\int_t^s G^n(r, \dbP_{Y_r})dr + V(s, {\dbP}_{Y_{s-}};\dbP_{\a^{t,l}_{s-}})\\
&&\dis= \int_t^s F(r,\hat\dbP_{\a^{t,l}_r}, \hat\dbP_{Y_r})dr-\sum_{\t^{\hat\dbP}_n \geq t}\int_t^s G^n(r, \hat\dbP_{Y_r})dr + J(s, \hat \dbP;\hat\dbP_{\a^{t,l}_{s-}})\\
&&\dis= J(t, \hat\dbP;l) \le V(t, m;l).
\end{eqnarray*}
Since $\dbP\in \cP(t,m;l)$ is arbitrary, we obtain $ \tilde V(t,m;l) \le  V(t,m;l)$. This completes the proof.
\qed

\subsection{Variational inequality}
We recall some differential calculus tools on the Wasserstein space. We say that a function $u:\cP_2(\bS) \to \dbR$ has a functional linear derivative $\delta_m u:\cP_2(\bS)\times\bS \to \dbR$ if
 \beaa
u(m')-u(m) 
= 
\int_0^1 \int_{\bS} \d_m u(\l m' + (1-\l)m, y)(m'-m)(dy)d\l
&\mbox{for all}&
m,m'\in\cP_2(\bS),
\eeaa
$\d_m u$ is continuous for the product topology, with $\cP_2(\bS)$ equipped with the 2-Wasserstein distance, and has quadratic growth in $x \in \dbR^d$, locally uniformly in $m \in \cP_2(\bS)$, so as to guarantee integrability in the last expression. Similarly, we define $\delta_{l,m} u:\cP_2(\bS \times \cM)\times\bS \times \cM \to \dbR$ for both $l$ and $m$ as follows
\beaa
u(m';l')-u(m;l) 
= 
\int_0^1 d\l\underset{q\in \cM}{\sum}\int_{\bS} \d_{l,m} u(\l m' + (1-\l)m, y;\l l' + (1-\l)l,q)(m'-m)(dy)(l'-l)(q).
\eeaa
 Denote $u^i(t,m;l)$ the part of $u(t,m;l)$ in the state $i$, we introduce the following measure flow generator.
\begin{equation}\label{DiffOp}
\left.\ba{c}
\dis\dbL u(t,m;l) := \sum_{q\in\cM}l(q) \sum_{i \in \cN}\Big\{\pa_t u^i(t,m;l) + \int_{\dbR^d} \cL^{i,q}_x\d_{l,m} u(t,m,x,i;l,q) m(dx,i)\Big\},\\
\dis \cL^{i,q}_x\d_{l,m} u(t,m,x,i;l):= b(q,x,i,m)\cd \pa_x \d_{l,m} u(t,m,x,i;l,q) + {1\over 2} \si^2(q,x,i,m): \pa_{xx}^2  \d_{l,m} u(t,m,x,i;l,q),\\
\dis \dbQ u(t,m;l) :=\sum_{p\in\cM}l(p)\int_{\bS}\sum_{q\neq p}\lambda_{pq}(x)\big[\d_{l,m} u(t,m,y;l,q)-\d_{l,m} u(t,m,y;l,p)\big]m(dy). 
\ea\right. 
\end{equation}

\begin{defn}
 $C_2^{1,2}([0, T]\times \cP_2(\dbR^{d'}))$ denotes the set of functions $u: [0, T]\times \cP_2(\dbR^{d'}) \to \dbR$ such that 

$\bullet$  $\pa_t u$, $\d_m u$, $\pa_y \d_m u$, $\pa_{yy}^2 \d_m u$ exist and are continuous in all variables; 

$\bullet$ $\pa_{yy}^2 \d_m u$ is bounded in $y$, locally uniformly in $(t, m)$.
\end{defn}

By abusing the notation, in the following statement, we let $Y$ denote a general c\`adl\`ag $\dbR^{d'}-$valued semimartingale on $[0, T]$. We denote $Y^c$  the continuous part of $Y$;  $Y^c_t = Y_0 + M^c_t + A^c_t$ the Doob-Meyer decomposition, where $M^c$ is the martingale part and $A^c$ is the finite variation part;   $\|A^c\|_t$  the total variation process of $A^c$ and $\la M^c\ra_t$ the quadratic variation process of $M^c$. We also denote $\a$ a general  $\cM-$valued Markov chain on $[0, T]$. Then we obtain the following It\^o's formula and relegate the proof to Appendix \ref{appendixA}.
\begin{prop}[It\^o's formula]\label{Itothm}
Let $u\in C^{1,2}_2(\overline{\mathbf{Q}}_0)$, and assume 
\bea\label{summability}
 \dbE\Big[ \|A^c\|_T^2 + \la M^c\ra_T  + \Big(\sum_{0 < s \le T} \lvert  Y_s-Y_{s-} \rvert \Big)^2 \Big] < \infty.
\eea
Then, denoting $\bm = \{m_s\}_{0\le s\le T}$  the marginal laws of $Y_s$, and  $\{l_s\}_{0\le s\le T}$ for $\a_s$
\begin{equation}
\label{JumpFunIto}
\left.\ba{lll}
\dis  u(T,m_{T-};l_{T-}) = u(0,m_{0-};l_{0-}) +  \int_0^T \Big[\pa_t u(s,m_s;l_s) + \dbQ u(s,m_s;l_s)\Big]ds\\ 
\dis\qq  + \dbE \int_0^T \big[ \pa_y \d_{l,m} u(s, m_s,Y_s;l_s,\a_s) \cd d A^c_s+  \pa_{yy}^2 \d_{l,m} u(s, m_s,Y_s;l_s,\a_s) : d\la Y^c\ra_s\big]\\
\dis\qq  + \sum_{s \in J_{[0, T)}(\mathbf{m})} \Big[u(s,m_s;l_s)-u(s, m_{s-};l_s)\Big] + \dbE\Big[ \sum_{s\in J_{(0, T]}^c(\bm)}  \!\!\!\!\big(\d_m u(s, m_s, Y_s;l_s)- \d_m u(s, m_s, Y_{s-};l_s)\big)\Big].
\ea\right.
\end{equation}
\end{prop} 
Here
\begin{equation}
\label{JAm}
J_{[0, T)}(\bm) := \{s \in [0, T): m_s \neq m_{s-}\},\q J^c_{[0, T)}(\bm) := \{s \in [0, T): m_s = m_{s-}\}.
\end{equation}
We emphasize that unlike $m$, the discontinuity part of $l$ is implicit in $\dbQ u(t,m;l)$ since the regime-switching of Markov chain. The difference is that $l$ transfers without control, while $m$ under control.
We next introduce a partial order $\preceq$ on $\cP_2(\mathbf{S})$: we say that $m' \preceq m$ if
\bea\label{order}
m'(dx, j) = \sum_{i \in \cN} p_{ij}(x) m(dx, i), \ \mbox{s.t.} \  \sum_{j \in \cN} p_{ij}(x)=1, 
\eea
for some measurable $p : \dbR^d \rightarrow [0,1]$. If there exists $i\in \cN$ s.t. $p_{ii}(x)\neq 1$, then we call that $m' \prec m$.Then the variational inequality corresponding to our mean field optimal switching problem is defined by
\bea\label{obstacle}
\dis  \min\Big\{-(\dbL u  + F + \dbQ u), u - M[u]\Big\}(t,m;l) = 0, 
\eea
with boundary condition $u|_{t=T} = H$. Here $M[u](t,m;l):=\underset{m' \prec m} {\sup}\Big\{u(t,m';l)-\underset{i \in \cN}{\sum}\underset{j \neq i}{\sum}\int_{\dbR^d}g_{ij}(t,x)p_{ij}m(dx,i)\Big\}$. The above variational inequality is an obstacle problem on Wasserstein space, and it's not difficult to observe that \reff{obstacle} deduces to the form of variational inequality in our previous work \cite{LW} without Markov chain. Moreover, \reff{obstacle} keeps consistent with the standard case  if there is no mean field term as we can see in the following example. 

We consider the case  that $b$, $\si$ do not depend on $\cP_2(\mathbf{S})$-valued variable, and define the optimal switching problem
\begin{eqnarray*}
V(t,m;l) := \sup_{\dbP \in \cP(t,m;l)} \dbE^\dbP[h(\a_T,X_T)+\int_t^T f(s,X_s,I_s,\a_s)ds-\sum_{\t^{\dbP}_n \geq t}g_{\xi_{n-1}\xi_n}(\t_n,X_{\t_n})],\q (t,m,l) \in \bQ_0.
\end{eqnarray*}
We  also introduce  $v^i(t,x;p):= V(t, \d_{(x,i)};\d_p)$  related to  the standard obstacle problem, see Tao, Wu  \& Zhang \cite{TWZ}.
\begin{equation}\label{varineq}
\left.\begin{array}{c}
\min\{-(\pa_t v^i +\cL^{i,p}_x v^i + f^i)-\underset{q \neq p}{\sum}\l_{pq}[v^i(t,x;q)-v^i(t,x;p)], v^i-\underset{j\neq i}{\max}(v^j-g_{ij})\}(t,x;p) = 0, \q v(T,\cd) = h,\\
\mbox{where}\q \cL^{i,p}_x v^i(t,x;p):=  b(x,i,p)\cd \pa_x v^i(t,x;p) + {1\over 2} (\si)^2(x,i,p): \pa_{xx}^2 v^i(t,x;p).
\end{array}\right.
\end{equation}

\begin{prop}\label{standardprop}
Assume $v\in C_2^{1,2}([0, T]\times \dbR^d)$. Then: 
 $V(t,m;l) = \underset{p\in \cM}{\sum}l(p)\underset{i\in \cN}{\sum}\int_{\dbR}v^i(t,x;p)m(dx,i)$, and $V$ is a classical solution of the corresponding obstacle equation on the Wasserstein space.
\end{prop}
\proof Obviously, $V \in C_2^{1,2}(\overline \bQ_0)$ with 
\begin{eqnarray*}
\pa_t V(t,m;l) = \underset{p\in \cM}{\sum}l(p)\underset{i\in \cN}{\sum}\int_{\dbR^d}\pa_t v^i(t,x;p)m(dx,i), \q \d_{l,m} V(t,m,x,i;l,p) = v^i(t,x;p). 
\end{eqnarray*}
We then show that $V$ is a solution of the equation \eqref{obstacle}.
First, terminal condition is satisfied naturally. Next, by \eqref{DiffOp} and \eqref{varineq}, we easily check
\begin{eqnarray*}
&(\dbL V + \dbQ V+F)(t,m;l)=\underset{p\in \cM}{\sum}l(p)\underset{i\in \cN}{\sum}\int_{\dbR^d}\big(\pa_t v^i +\cL^{i,p}_x v^i + f^i+\underset{q \neq p}{\sum}\l_{pq}(x)[v^i(t,x;q)-v^i(t,x;p)]\big)m(dx,i).
\end{eqnarray*}
This proves the first part of variational inequality.
 Define $A^{i,p}_t := \{x : v^i(t,x;p) > v^j(t,x;p)- g_{ij}(t,x) \}$ and $m^{A^{i,p}_t} := m \circ (\mathbf{x}, i \1_{A^{i,p}_t}(\mathbf{x})+ j \1_{(A^{i,p}_t)^c}(\mathbf{x}))^{-1}$, then for $m=\d_{(x,i)}$, $l=\d_{p}$
\begin{equation}
\begin{split}
& V(t,m;l) - V(t,m^{A^{i,p}_t};l)\\
& = \underset{p\in \cM}{\sum}l(p)\big\{\int_{\dbR^d}v^i(t,x;p)m(dx,i)-\int_{A_t^{i,p}}v^j(t,x;p) m(dx,i)-\int_{(A_t^{i,p})^c}v^j(t,x;p) m(dx,j)\big\}\\
&\ge \underset{p\in \cM}{\sum}l(p)\big\{-\int_{A_t^{i,p}}g_{ij}(t,x) m(dx,i)-\int_{(A_t^{i,p})^c}g_{ij}(t,x) m(dx,j)\big\}=-\underset{p\in \cM}{\sum}l(p)\int_{\dbR^d}g_{ij}(t,x)m^{A^{i,p}_t}(dx), 
\end{split}
\end{equation}
 $"="$ holds if and only if $v^i(t,x;p) = v^j(t,x;p)- g_{ij}(t,x)$, and $p_{ij}=1$ at this time. This proves the second part of variational inequality. 
Thus, $V$ is a solution of \eqref{obstacle} since $v$ is the solution of \eqref{varineq}.
\qed
\section{Viscosity solutions }\label{sect-viscosity}

\subsection{Viscosity property}

For $\d > 0$ and $(t,m,l) \in \mathbf{Q}_0$, we introduce the neighborhood
$$ \cN_\d (t,m;l) := \big\{(s, \tilde m) :  s \in [t, t+\d], \dbP \in \cP(t,m;l), \ \tilde m \in \{ \dbP_{Y_{s-}}, \dbP_{Y_s} \} \big\}. $$
Note that, as the closure of a c\`adl\`ag $\cP_2(\bS)$-valued graph, $\cN_\d(t,m;l)$ is compact, by the compactnesses of $[t,t+\d]$, $\cP(t,m;l)$ and $\{ (\dbP_{Y_{s-}}, \dbP_{Y_s}) \}_{s \in [t,t+\d]}$ for any $\dbP \in \cP(t,m;l)$. 

For a  locally bounded function $u: \overline{\mathbf{Q}}_0\longrightarrow\dbR$, we introduce its LSC and USC envelopes relatively to $\cP(t,m;l)$, $u_*$ and $u^*$ respectively:
\beaa
u_*(t,m;l) := \underset{(s, \tilde m) \to (t,m)}{\lim \inf} u(s,\tilde m;l), \q u^*(t, m;l) := \underset{(s,\tilde m) \to (t,m)}{\lim \sup} u(s,\tilde m;l), 
&\mbox{for all}&
(t,m,l) \in \overline{\mathbf{Q}}_0, 
\eeaa
where the limits are taken on all sequences $\{t_n,m_n\}_{n \ge 1}$ converging to $(t,m)$, with $(t_n,m_n) \in \cN_{T-t}(t, m;l)$ 
for all $n$. We then introduce the sets of test functions
\bea
\overline{\cA}u(t,m;l) &:=& \Big\{\f \in C^{1,2}_2([t,T]\times\cP_2(\mathbf{S})):  (\f-u_*)(t,m;l) = \max_{\cN_\d(t, m;l)}(\f-u_*)(s,\tilde m;l)\Big\} \nonumber, \\ 
\underline{\cA}u(t,m;l) &:=& \Big\{\f \in C^{1,2}_2([t,T]\times\cP_2(\mathbf{S})): (\f-u^*)(t,m;l) = \min_{\cN_\d(t, m;l)}(\f-u^*)(s,\tilde m;l) \Big\}. \nonumber
\eea
 It should be noted that $\f$ is a function of $(t,m)$ rather than $l$ here. As we can see in the definition below, $\f$ is only a smooth replacement of $\dbL u(\cd,\cd;l)$ for given $l$.

\begin{defn}
\label{defn-viscosity}
Let $u : \mathbf{Q}_0 \rightarrow \dbR$ be locally bounded.

\no{\rm (i)} $u$ is a viscosity supersolution of \reff{obstacle} if, for any  $(t, m,l)\in\mathbf{Q}_0$,
\bea
\label{super}
\Big(u_* - M[u_*]\Big) (t,m;l) \ge 0,\q
\mbox{and}\q
-(\dbL \f + F + \dbQ u_* ) (t,m;l) \ge 0,
~\forall \f \in \overline{\cA}u(t,m;l). 
 \eea
 
 \no{\rm (ii)} $u$ is a viscosity subsolution of \reff{obstacle} if,  for any  $(t, m,l)\in\mathbf{Q}_0$ ,
\bea
\label{sub}
\min\Big\{-(\dbL \f  + F + \dbQ u^*), u^* - M[u^*]\Big\}(t,m;l) \le 0,\q\forall~
 \f \in \underline{\cA}u(t,m;l).
 \eea

\no{\rm (iii)} $u$ is a viscosity solution of \eqref{obstacle} if it is a viscosity supersolution and subsolution.

\end{defn}
\begin{rem}
 Actually, there are two ways of definition by taking the neighborhood either in $(t,m)$ or $(t,m,l)$. The convenience of the former is under which we can apply the two-time-scale approach in \S\ref{sect-TTS}. 
\end{rem}

\begin{thm}
\label{thm-existence}
The value function $V$  is a viscosity solution of \eqref{obstacle}.
\end{thm}
\proof First, under Assumption \ref{assum-bsig},  $V$ inherits the local boundedness of $h$.

 \no(i) We first verify the viscosity supersolution property. Fix $(t, m,l)$ and $\f\in \overline{\cA}V(t,m;l)$. We may assume w.l.o.g. that $[V_*-\f](t,m;l) = 0$. Let $\d > 0$ and $(t_n, m_n)_{n \ge 1} \in \cN_\d(t,m;l)$ converging to $(t,m)$ s.t. $V(t_n, m_n;l) \underset{n \to \infty}{\longrightarrow} V_*(t,m;l)$, and denote $\eta_n := [V-\f](t_n,m_n;l) \ge 0$, as $V\ge V_*$. Thus, we have $\eta_n \underset{n \to \infty}{\longrightarrow} 0$. By the DPP \eqref{weakDPP}, we have 
\bea\label{exist-dpp}
\eta_n + \f(t_n, m_n)
&=&
V(t_n, m_n;l) 
\ge \int_{t_n}^{s_n} F(r, \bar \dbP^{m_n;l}_{\a^{t_n,l}_r},\bar \dbP_{Y_r}^{m_n;l})dr + V(s_n, \bar \dbP_{Y_{s_n-}}^{m_n;l};\bar \dbP_{\a^{t_n,l}_{s_n-}}^{m_n;l})  \nonumber\\
&\ge& \int_{t_n}^{s_n} F(r, \bar \dbP^{m_n;l}_{\a^{t_n,l}_r},\bar \dbP_{Y_r}^{m_n;l})dr + V_*(s_n, \bar \dbP_{Y_{s_n-}}^{m_n;l};\bar \dbP^{m_n;l}_{\a^{t_n,l}_{s_n-}}), 
\eea
where $\bar \dbP^{m_n;l}:=\bar \dbP^{t_n, m_n;l} \in \cP(t_n, m_n;l)$ is defined by \reff{barP} such that $X$ is unstopped, and $s_n := t_n + h_n$ with $h_n := \sqrt{\eta_n} \vee n^{-1}$. Define
\begin{equation}
\psi(s,\tilde m;\tilde l)=\f(s,\tilde m)+V_*(s_n, \bar \dbP_{Y_{s_n-}}^{m_n;l};\tilde l)-\f(t_n, m_n),
\end{equation}
then by It\^o's formula, the above gives
\bea\label{exist-ito}
&&\dis \psi(s_n, \bar \dbP_{Y_{s_n-}}^{m_n;l};\bar \dbP^{m_n;l}_{\a^{t_n,l}_{s_n-}})=\psi(t_n, m_n;l)+\int_{t_n}^{s_n} (\dbL \psi  + \dbQ \psi)(r, \bar \dbP_{Y_{r}}^{m_n;l};\bar \dbP^{m_n;l}_{\a^{t_n,l}_r})dr\nonumber\\
&&\dis =V_*(s_n, \bar \dbP_{Y_{s_n-}}^{m_n;l};l)+ \int_{t_n}^{s_n} \big(\dbL \f(r, \bar \dbP_{Y_{r}}^{m_n;l}) + \dbQ V_*(s_n, \bar \dbP_{Y_{s_n-}}^{m_n;l};\bar \dbP^{m_n;l}_{\a^{t_n,l}_r}) \big)dr.\\
&&\dis \geq \f(s_n, \bar \dbP_{Y_{s_n-}}^{m_n;l})+ \int_{t_n}^{s_n} \big(\dbL \f(r, \bar \dbP_{Y_{r}}^{m_n;l}) + \dbQ V_*(s_n, \bar \dbP_{Y_{s_n-}}^{m_n;l};\bar \dbP^{m_n;l}_{\a^{t_n,l}_r}) \big)dr \nonumber.
\eea
Then by \reff{exist-dpp} and \reff{exist-ito} we have
\beaa
\eta_n \geq \int_{t_n}^{s_n} \big( F(r, \bar \dbP^{m_n;l}_{\a^{t_n,l}_r},\bar \dbP_{Y_r}^{m_n;l})+\dbL \f(r, \bar \dbP_{Y_{r}}^{m_n;l})+\dbQ V_*(s_n, \bar \dbP_{Y_{s_n-}}^{m_n;l};\bar \dbP^{m_n;l}_{\a^{t_n,l}_r})\big)dr.
\eeaa 
Send $n\to \infty$, since $h_n{\longrightarrow} 0$, we obtain  $-(\dbL \f + F + \dbQ V_*)(t,m;l) \ge 0$. 

We now prove the remaining part of the supersolution property. Let $m' \preceq m$ with transition probability $p$, we may assume without loss of generality that $p$ is continuous. For all $n \ge 1$, define $m_n' \preceq m_n$ as the measure obtained from $m_n$ by applying the same $p$. Given the continuity of $p$ and the compactness $\cN_\d(t,m;l)$, we see by \eqref{order} that $\cW_2(m_n', m') \underset{n \to \infty}{\longrightarrow} 0$. Let $\dbP^{m_n, m_n';l} \in \cP(t_n, m_n;l)$ be s.t. $ \dbP_{Y_{t_n}}^{m_n, m_n';l} = m_n'$, 
 \bea\label{superexistence}
 V(t_n, m_n;l) \ge  V(t_n, m_n';l)-\underset{i \in \cN}{\sum}\underset{j \in \cN}{\sum}\int_{\dbR^d}g_{ij}(t,x)p_{ij}m_n(dx,i).
 \eea
Take $\lim \inf_{n\to\infty}$ in \reff{superexistence}, we obtain $ V_*(t,m;l) \ge M[V_*](t,m;l)$ as $V(t_n, m_n;l) \to V_*(t,m;l)$. 

\no(ii) We next verify the viscosity subsolution property.   Let  $\f \in \ul \cA V(t,m;l)$. We may assume w.l.o.g. that $[V^*-\f](t,m;l) = 0$.  Let  $(t_n, m_n) \in \cN_\d(t,m;l)$ converging to $(t,m)$ s.t. $V(t_n, m_n;l) \underset{n \to \infty}{\longrightarrow} V^*(t,m;l)$, and denote $-\eta_n := [V-\f](t_n,m_n;l) \le 0$.  Thus $\eta_n \underset{n \to \infty}{\longrightarrow} 0$. Suppose $V^* (t,m;l)>M[V^*](t,m;l)$, then $V(t_n,m_n;l)>M[V](t_n,m_n;l)$  for $n$ large enough. Since $\cP(t_n, m_n;l)$ is compact, there exists  $\dbP^{n,*}$ and $s_n > t_n$ s.t. $V(t_n, m_n;l) = \int_{t_n}^{s_n} F(r,\dbP_{\a^{t_n,l}_r}^{n,*},\dbP_{Y_r}^{n,*})dr + V(s_n,\dbP_{Y_{s_n-}}^{n,*};\dbP_{\a^{t_n,l}_{s_n-}}^{n,*})$. Then same with (i), we have
\bea\label{subexistence}
 -\eta_n \le \int_{t_n}^{s_n} \big(F(r,\dbP_{\a^{t_n,l}_r}^{n,*}, \dbP_{Y_r}^{n,*})+\dbL \f(r, \dbP_{Y_{r}}^{n,*})+\dbQ V^*(s_n,\dbP_{Y_{s_n-}}^{n,*};\dbP_{\a^{t_n,l}_r}^{n,*})\big)dr. 
 \eea
 For the trajectories are c\`adl\`ag, we obtain $-(\dbL \f + F + \dbQ V^*) (t,m;l) \le 0$ by sending $n\to \infty$. \qed

\subsection{Regularity and stability}
We first state some regularity results which will be used in the rest of this section. First up, we state the greater continuity in the following lemma, the proof of which is postponed to the Appendix \ref{appendixB}.

\begin{lem}
\label{thm-reg}
Assume $f, h, G^n$ are uniformly continuous in $(t, y, m)$ and $b, \si$ are uniformly Lipschitz continuous in $m$ under $\cW_1$,   then $V$ is also continuous in $m$ under $\cW_1$.  
\end{lem}

We note that $V(t,m;l)$ is naturally continuous in $l$ under $\cW_1$, because state process $Y$ is not affected by initial $l$, but evolves likewise on each $\a\in \cM$.  Therefore, $\dbP_Y$ is continuous in $l$ under distribution distance.
For further smoothness requirements, one can hardly expect it to be true even for the standard optimal switching problems. We next establish a regularity result for the value function without switching.  For $(t,m,l) \in \bQ_0$, let $\bar \dbP^{t, m} \in \cP(t,m;l)$ be as by \reff{barP}, and define
\bea
\label{XU}
 U(t, m;l) := H(\bar\dbP^{t,m}_{\a^{t,l}_T},\bar \dbP_{Y_T}^{t,m})  + \int_t^T F(r, \bar\dbP^{t,m}_{\a^{t,l}_r},\bar\dbP_{Y_r}^{t,m})dr.
\eea
It's easy to check the continuity of $U$, since there is no switching.
For $\f=b, \si$,  assume  $\pa_x\f, \d_m \f, \pa_x \d_m \f, \pa_{xx}^2\d_m\f$ exist and are continuous and bounded  and  all ther derivatives of $\f$ are Lipschitz up to order 2.  Then $U \in C^{1,2}(\overline\bQ_0)$ with bounded  $\pa_x \d_m U, \pa_{xx}\d_m U$. Moreover, if $b, \si$ are  uniformly Lipschitz continuous in $m$ under $\cW_1$, then $U$ is uniformly Lipschitz continuous in $m$ under $\cW_1$.

Given the continuity of $U$, we introduce a smooth mollifier for functions on the Wasserstein space. 

  \begin{lem}\label{lem-mol}
{\rm (i)} Let $U : \cP_2(\bS) \rightarrow \dbR$ be continuous. There exists $\{U_n\}_{n \ge 1}$ in $C^\infty(\cP_2(\bS))$ such that $\dis\lim_{n\to\infty} \sup_{m\in \cK}|U_n(m) - U(m)| =0$ for any compact set  $\cK\subset \cP_2(\bS)$;\\
  {\rm (ii)} Let $U : \cP_1(\bS) \rightarrow \dbR$ be continuous under $\cW_1$. There exists $\{U_n\}_{n \ge 1}$ in $C^\infty(\cP_2(\bS))\cap C^0(\cP_1(\bS))$ such that  $\dis\lim_{n\to\infty} \sup_{m\in \cK}|U_n(m) - U(m)| =0$ for any compact  set $\cK\subset \cP_1(\bS)$;\\
{\rm (iii)} Assume further that $U$ is Lipschitz continuous under $\cW_1$, then we may choose $\{U_n\}_{n \ge 1}$ to be Lipschitz continuous under $\cW_1$, uniformly in $n$. 
\end{lem}

The mollifier is adopted from Mou \& Zhang \cite{MZ}. Note that the extension of the state space from $\dbR^d$ in \cite{MZ} to $\bS$ here is straightforward, so we omit the proof. Now, we state the main result of this subsection.
\begin{thm}
Let  $\{u_\e\}_{\e > 0}$ and $\{v_\e\}_{\e > 0}$ be two families of viscosity subsolutions and supersolutions of \eqref{obstacle}, respectively. Assume that the following relaxed semilimits are finite
\beaa
\ol u(t,m;l) := \underset{(s,\tilde m) \rightarrow (t,m)}{\lim \sup} u_\e(s,\tilde m;l),
~\mbox{and}~
\ul v(t,m;l) := \underset{(s,\tilde m) \rightarrow (t,m)}{\lim \inf} v_\e(s,\tilde m;l),
&(t,m,l) \in \overline{\mathbf{Q}}_0,&
\eeaa
where the limits are sequences $(\e_n, t_n, m_n) \rightarrow (0,t,m)$, with $(t_n,m_n) \in \cN_{T-t}(t,m;l)$. Then $\ol u$ (resp. $\ul v$) is a USC (resp. LSC) viscosity subsolution (resp. supersolution) of \eqref{obstacle}.  
\end{thm}
\proof
(i) We prove the stability of the supersolution first. Observe that we may assume without loss of generality that $v_\eps$ is LSC as $
\ul v(t,m;l) = \underset{(s,\tilde m) \rightarrow (t,m)}{\lim \inf} (v_\e)_*(s,\tilde m;l)$.

Fix $(t,m,l) \in {\mathbf{Q}}_0$, and $\f \in \ol \cA \ul v(t,m;l)$ with corresponding $\d$, and s.t. $(t,m)$ is a strict maximizer of $\f - \ul v$ on $\cN_\d(t,m;l)$. By definition, there exists a sequence $(t_n, m_n)\to (t, m)$ s.t. $v_{\e_n}(t_n,m_n;l) \rightarrow \ul v(t,m;l)$. Note that $(t_n, m_n) \in \cN_\d(t,m;l)$ for all $n$ large, then we can find $\d' < \d$ s.t. $\cN_{\d'}(t_n,m_n;l) \subset \cN_\d(t,m;l)$. Let $(\hat t_n, \hat m_n)$ be a maximizer of $\f -v_{\e_n}$ on $\cN_{\d'}(t_n,m_n;l)$.  We first note that
\bea\label{conv-stable}
(\hat t_n, \hat m_n) \underset{n \rightarrow \infty}{\longrightarrow} (t,m).
\eea
Indeed,  $(\hat t_n, \hat m_n) \in \cN_{\d'}(t_n,m_n;l) \subset \cN_\d(t,m;l)$ for all $n$. Thus, by compactness, there exists a subsequence (still named $\hat m_n$) converging to some $(\hat t, \hat m) \in \cN_\d(t,m;l)$. Observing that 
\beaa
[\f - \ul v](t,m;l) &=& \underset{n \rightarrow \infty}{\lim}[\f - v_{\e_n}](t_n,m_n;l) \le \underset{n \rightarrow \infty}{\lim \inf}[\f - v_{\e_n}](\hat t_n, \hat m_n;l)\\
&\le&  \underset{n \rightarrow \infty}{\lim \sup}[\f - v_{\e_n}](\hat t_n, \hat m_n;l) \le [\f - \ul v](\hat t, \hat m;l),
\eeaa
we conclude by $(t,m)$ is a strict maximizer of $\f - \ul v$ on $\cN_\d(t,m;l)$ that $(\hat t, \hat m) = (t,m)$, and thus \eqref{conv-stable} holds true. Given that $(t_n, m_n)$ and $(\hat t_n, \hat m_n)$ have the same limit, we have $\cN_{\d''}(\hat t_n, \hat m_n;l) \subset \cN_{\d'}(t_n, m_n;l)$ for some $\d'' < \d'$ and $n$ large enough. As $(\hat t_n, \hat m_n)$ is also a maximizer on $\cN_{\d''}(\hat t_n, \hat m_n;l)$, the supersolution property implies
$ -(\dbL \f + F + \dbQ v_{\e_n}) (\hat t_n, \hat m_n;l) \ge 0,$  
and we derive the first part of the supersolution property of $\ul v$ by sending $n \rightarrow \infty$. 

We now prove that $\big(\ul v - M[\ul v]\big)(t,m;l) \ge 0$. By the supersolution property, we have 
$\big(v_{\e_n} - M[v_{\e_n}]\big) (t_n, m_n;l) \ge 0$, for all $n \ge 1$,
 and we conclude by taking the $\lim \inf$ that $\big( \ul v - M[\ul v]\big)(t,m;l) \ge 0$, as the l.h.s. of the inequality converges. 
 
\no(ii) We now prove the stability of the subsolution. Similarly to (i), we may assume that $\{u_\e\}_{\e > 0}$ is a family of USC viscosity subsolutions of \eqref{obstacle}.
 Let $(t,m,l)\in\bQ_0$ and $\f \in \ul \cA \ol u(t,m;l)$ be such that $(t,m)$ is a strict local minimizer of $\varphi-\overline{u}$.  Following the same argument as in the previous step, replacing maximizers with minimizers, we may construct $(\hat t_n, \hat m_n)$,  converging to some $(\hat t,\hat m)$, and satisfying the inequalities
$$[\f - \ol u](t,m;l)  \ge \underset{n \rightarrow \infty}{\lim \sup} \ [\f - u_{\e_n}](\hat t_n, \hat m_n;l) \ge \underset{n \rightarrow \infty}{\lim \inf} \ [\f - u_{\e_n}](\hat t_n, \hat m_n;l) \ge [\f - \ol u](\hat t, \hat m;l). $$
By the strict minimum property of $(t,m)$, this again implies that $(\hat t,\hat m)=(t,m)$, and $\underset{n \rightarrow \infty}{\lim} u_{\e_n}(\hat t_n, \hat m_n;l) = \ol u(t,m;l)$.
Assume $\big( \ol u - M[\ol u]\big)(t,m;l)>0$, then $\big(u_{\e_n} - M[u_{\e_n}]\big) (\hat t_n, \hat m_n;l)>0$ for $n$ enough large.
 Then for $\f \in \ul \cA u_{\e_n}(\hat t_n,\hat m_n;l)$, the viscosity subsolution implies $-(\dbL \f + F + \dbQ u_{\e_n}) (\hat t_n, \hat m_n;l) \le 0$ for $n$ large enough, and we conclude by letting $n \longrightarrow \infty$. \qed

\subsection{Comparison}
The uniqueness is a direct corollary of the following comparison. Thus, below is the main results of this subsection.

\begin{thm}[Comparison]
\label{thm-comparison}
Assume that $b$, $\si$, $f$, $h$ can be extended to $\cP_1(\bS)$ under $\cW_1$ continuously; $b$, $\si$ are uniformly Lipschitz continuous in $(x, m)$ under $\cW_1$. 

{\rm (i)}  Let $v$ be a LSC viscosity supersolution of \reff{obstacle} satisfying $v|_{t=T} \ge H$.  Then $v \ge V$. 

{\rm (ii)} Let $u$ be an USC viscosity subsolution of \reff{obstacle} satisfying $u|_{t=T} \le H$. Then $u\le V$.
\end{thm}
The proof is lengthy, let's start with a brief summary of the idea. We assume the comparison does not hold, then we can construct some  test functions by $U$ defined in \reff{XU} and  derive a contradiction  to the nature of the viscosity solutions. Thanks to mollifier in Lemma \ref{lem-mol}, it is possible to ensure the test functions smooth enough .

\proof (i) We first compare $V$ and $v$. Fix $\e>0$. For each $n\ge 1$, denote $t_i := t^{(n)}_i:=\frac {iT} n$, $0 \le i \le n$. First, note that, for $(t,m,l) \in \mathbf{Q}_0$, it follows from the continuity of the coefficients that
\small{
\bea
\label{Vn}
\left.\ba{c}
\dis V(t, m;l) := \lim_{n\to\infty} V_n(t, m;l)= \lim_{n\to\infty}\sup_{\dbP\in \cP_n(t, m;l)} \Big\{ \int_t^T F(r,\dbP_{\a^{t,l}_r},\dbP_{Y_r})dr + H(\dbP_{\a^{t,l}_T},\dbP_{Y_T})-\sum_{\t^{\dbP}_k \geq t}\int_t^T G^k(r, \dbP_{Y_r})dr \Big\},\\
\dis \cP_n(t, m;l):= \Big\{\dbP\in \cP(t, m;l): \mbox{switching time $\t^\dbP$ takes values in $\{t_1,\ldots,t_n\} \cap [t, T]$, $\dbP$-a.s.}\Big\}.
\ea\right. \nonumber
\eea
}

{\it Step 1:} We show that $(V_n - v)(t_{n-1}, \cdot;\cdot) \le \frac \e n$. Assume to the contrary that there exists $(m_{n-1},l_{n-1})$ such that $(V_n-v)(t_{n-1}, m_{n-1};l_{n-1}) > \frac \e n$. By the definition of $\cP_n(t_{n-1}, m_{n-1};l_{n-1})$, for subsequent $(m,l)$ we have $V_{n}(t, m;l)= \int_t^T F(r, \bar \dbP^{t,m}_{\a^{t,l}_r}, \bar \dbP^{t,m}_{Y_r})dr + H(\bar \dbP^{t,m}_{\a^{t,l}_T},\bar \dbP_{Y_T}^{t,m})$, $t\in (t_{n-1}, T]$,  where $\bar \dbP^{t, m} \in \cP(t,m;l)$ is defined by \reff{barP}. 

Let $\d_1, \d_2>0$ be small numbers which will be specified later. Applying Lemma \ref{lem-mol} (ii), (iii), let $(h_k, f_k, b_k)$ be the smooth mollifier of $(g, f, b)$ (under $\cW_1$), where $b_k$ is also mollified in $(t,x)$ in a standard way, such that $\lVert h_k-h \rVert_\infty + \lVert f_k-f \rVert_\infty  \le \d_1$, $\lVert b_k - b \rVert_\infty \le \d_2$, and $g_k$ is Lipschitz continuous under $\cW_1$ with a Lipschitz constant $L_k$ depending on $k$, and $b_k$ is uniformly Lipschitz continuous in $(x, m)$ under $\cW_1$ with a Lipschitz constant $L$ independent of $k$. By otherwise choosing a larger $L$ we assume $\si$ is also uniformly Lipschitz continuous in $(x, m)$ under $\cW_1$ with Lipschitz constant $L$. Let $U^{k_1, k_2}$ be defined by \reff{XU} corresponding to $(b_{k_2}, \si, h_{k_1}, f_{k_1})$. Then we have, 
\begin{equation}\label{eqdiff}
\pa_t U^{k_1,k_2}(t,m;l) + \dbQ U^{k_1,k_2}(t,m;l) + \underset{q\in\cM}{\sum}l(q)\int_{\bS} \Big[ b_{k_2}\cd \pa_x \d_m U^{k_1,k_2} + {1\over 2} \si^2 : \pa_{xx}^2 \d_m U^{k_1,k_2} + f_{k_1} \Big] m(dy) =0,
\end{equation}
and $U^{k_1,k_2}$ is Lipschitz continuous in $m$ under $\cW_1$ with a Lipschitz constant $C_{L, L_{k_1}}$ independent of $k_2$. This, in particular, implies $|\pa_x \d_m U_1^{k_1,k_2}|\le C_{L, L_{k_1}}$ for all $k_2$. Then, we deduce from \eqref{eqdiff} that
\bea
\label{Uk2est}
\Big| (\dbL U^{k_1,k_2} + F_{k_1} + \dbQ U^{k_1,k_2})(t, m;l) \Big| = \Big| \underset{q \in \cM}{\sum}l(q)\int_{\bS} (b - b_{k_2}) \cdot \pa_x \d_m U^{k_1,k_2} m(dy) \Big|   \le C_{L, L_{k_1}} \d_2,~\forall k_2 \ge 1,
\eea
where $F_{k_1}(t,m,l) := \underset{q \in \cM}{\sum}l(q)\int_{\bS} f_{k_1}(t,q,y,m)m(dy)$ as in \reff{F}. Moreover, since
\beaa
 U^{k_1, k_2}(t,m;l) =  H_{k_1}( \bar \dbP_{\a^{t,l}_T}^{t,m,k_2},\bar \dbP_{Y_T}^{t,m,k_2}) + \int_t^T F_{k_1}(r, \bar \dbP_{\a^{t,l}_r}^{t,m,k_2}, \bar \dbP_{Y_r}^{t,m,k_2}) dr,
\eeaa
where $\bar \dbP^{t,m,k_2}$ is s.t. $X$ is unstopped with drift coefficient $b_{k_2}$ instead of $b$, one can easily show that
\beaa
\Big|U^{k_1,k_2}(t,m;l) -\Big[H(\bar \dbP_{\a^{t,l}_T}^{t,m},\bar \dbP_{Y_T}^{t,m}) + \int_t^T F(r, \bar \dbP_{\a^{t,l}_r}^{t,m}, \bar \dbP_{Y_r}^{t,m}) dr \Big]\Big|\le C[\d_1 + \d_2] \le {\e\over 4n},
\eeaa
for $\d_1, \d_2$ small enough. Then
\beaa
&& V_n(t_{n-1}, m_{n-1};l_{n-1}) \\
&&= \sup_{m' \preceq m_{n-1}} \Big\{ \int_{t_{n-1}}^T F(r,\bar \dbP_{\a^{t_{n-1},l_{n-1}}_r}^{t_{n-1},m'}, \bar \dbP_{Y_r}^{t_{n-1},m'}) dr + H(\bar \dbP_{\a^{t_{n-1},l_{n-1}}_T}^{t_{n-1},m'},\bar \dbP_{Y_T}^{t_{n-1},m'} )-\underset{i \in \cN}{\sum}\underset{j \in \cN}{\sum}\int_{\dbR^d} g_{ij}p_{ij}m_{n-1}(dx,i) \Big\} \\
&&\le \sup_{m' \preceq m_{n-1}} U^{k_1, k_2}(t_{n-1}, m';l_{n-1}) -\underset{i \in \cN}{\sum}\underset{j \in \cN}{\sum}\int_{\dbR^d} g_{ij}p_{ij}m_{n-1}(dx,i) +{\e\over 4n}. 
\eeaa
By the supersolution property, $ \big(v-M[v]\big)(t,m;l) \ge 0$ implies 
$$v(t_{n-1}, m_{n-1};l_{n-1}) \ge v(t_{n-1}, m';l_{n-1})-\underset{i \in \cN}{\sum}\underset{j \in \cN}{\sum}\int_{\dbR^d} g_{ij}p_{ij}m_{n-1}(dx,i),$$ 
for $n$ large enough, hence 
\beaa
\frac \e n \le (V_n  - v)(t_{n-1}, m_{n-1};l_{n-1})  \le \sup_{m'\preceq m_{n-1}} (U_{k_1,k_2}-v)(t_{n-1}, m';l_{n-1}) +{\e\over 4n}. 
\eeaa
This implies that, 
\bea\label{supersol-test}
\max_{(s, \tilde m)\in \cN_{\frac T n}(t_{n-1}, m_{n-1};l_{n-1})} \Big\{(U^{k_1,k_2} - v)(s, \tilde m;l_{n-1}) - {T-s\over n}  \Big\} \ge \frac{3\e}{4n} - {T \over n^2}  \ge {\e \over 2n},  
\eea
for $n$ sufficiently large. Note that $(U^{k_1,k_2} - v)(\dbP_{\a_T}, \dbP_{Y_T}) \le (H_{k_1} - H)(\dbP_{\a_T},\dbP_{Y_T}) \le {\e \over 4n}$ for all $\dbP \in \cP(t_{n-1}, m_{n-1};l_{n-1})$ and $v$ is LSC, then by compactness of $\cN_{\frac T n}(t_{n-1}, m_{n-1};l_{n-1})$ there exists an optimal argument $(t^*, m^*), \ t^* < T$, to the above maximum. Thus $\f(s, \tilde m) := U^{k_1,k_2}(s,\tilde m;l_{n-1}) - {T-s\over n} \in \ol \cA v(t^*,m^*;l_{n-1})$, and therefore,
\beaa
0 &\le& -(\dbL \f + F + \dbQ v)(t^*,m^*;l_{n-1})\\*
 &=& \Big\{-(\dbL U^{k_1,k_2} + F_{k_1} + \dbQ U^{k_1,k_2})  + (F_{k_1} - F) + (\dbQ U^{k_1,k_2}- \dbQ v)\Big\}(t^*,m^*;l_{n-1}) - {1\over n} \\*
&\le& C_{L, L_{k_1}}\d_2 +  (F_{k_1} - F)(t^*,m^*,l_{n-1})  + (\dbQ U^{k_1,k_2}- \dbQ v)(t^*,m^*;l_{n-1}) - {1\over n},
\eeaa
where the last inequality thanks to \reff{Uk2est}.
Fix $k_1$  so that $(F_{k_1} - F)(t^*,m^*,l_{n-1}) \le \frac{1}{2n}$ and set $\d_2$ small enough, we obtain the desired contradiction. 
\ms

{\it Step 2:} We show that $(V_n - v)(t_{n-2}, \cdot;\cdot) \le {2\e \over n}$. Assume to the contrary that there exists $(m_{n-2},l_{n-2})$ such that $(V_n-v)(t_{n-2}, m_{n-2};l_{n-2}) > {2\e \over n}$.  By the DPP, we have
$$V_n(t_{n-2}, m_{n-2};l_{n-2}) = \sup_{\dbP \in \cP_n(t_{n-2}, m_{n-2};l_{n-2})} \dbE \Big\{ \int_{t_{n-2}}^{t_{n-1}} F(r,\dbP_{\a_r}, \dbP_{Y_r})dr + V_n(t_{n-1}, \dbP_{Y_{(t_{n-1})-}};\dbP_{\a_{(t_{n-1})-}})\Big\},$$
Observe the fact that $v$ is a viscosity supersolution of \eqref{obstacle} also implies that $v + \frac \e n$ is a viscosity supersolution. Moreover, by {\it Step 1}, we have $(v + \frac \e n)(t_{n-1}, \cdot;\cdot) \ge V_n(t_{n-1}, \cdot;\cdot)$. Thus, using the same procedure as {\it Step 1} (where $V_n$ replaces $H$ on $(t_{n-2}, t_{n-1}]$), it follows that 
$$ \big(V_n - (v+\frac \e n)\big)(t_{n-2}, \cdot;\cdot) \le \frac \e n. $$ 
Finally, by backwards induction, we have $(V_n - v)(t_{n-j}, \cdot;\cdot) \le {j\e \over n}$ for all $j\in \{0, \dots, n\}$, and thus
$ (V_n - v)(t, \cdot;\cdot) \le \e, $
which implies by the arbitrariness of $n$ and $\e$ that $v \ge V$.

\medskip

\noindent (ii)  We next compare $V$ and $u$. By the subsolution property of $u$, we divide the discussion into the following two cases. One is $ \big(u-M[u]\big)(t,m;l)\le 0$, and $-(\dbL \f  + F + \dbQ u)(t,m;l)> 0$, then we can prove by following the same procedure in (i). The other case is 
\bea\label{u-contra}
-(\dbL \f  + F + \dbQ u)(t,m;l)\le 0,\q \forall \f \in \ul \cA v(t,m;l).
\eea
 For this case, assume by contradiction that  
$
u(t, m;l) > V(t,m;l) 
$ for some $(t, m,l)$. Then, for $\e>0$ small enough, 
\bea
\label{utTe}
u(t, m;l) - \f_\e(t,m) > \sup_{\dbP\in \cP(t, m;l)} \Big\{ \int_t^T F(r,\dbP_{\a^{t,l}_r}, \dbP_{Y_r})dr +  H(\dbP_{\a^{t,l}_T},\dbP_{Y_T})- \sum_{\t^{\dbP}_k \geq t}\int_t^T G^k(r, \dbP_{Y_r})dr \Big\}, 
\eea
where $\f_\e(s, \tilde m) := \e \big[(T-s)+ \sum_{i \in \cN} \tilde m(\dbR^d,i)\big]$. Let $(t^*, m^*,l^*)$ be s.t.
\begin{equation}
\label{c*}
\dis (u-\f_\e)(t^*,m^*;l^*) + \int_t^{t^*} \!\!\! F(r, \bar\dbP_{\a^{t,l}_r},\bar\dbP_{Y_r})dr 
\dis = \!\!\! \max_{\tiny \begin{array}{c} (s, \tilde m)\in \cN_{T-t}(t,m;l)  \end{array}} \!\!\! \Big\{(u-\f_\e)(s, \tilde m;\bar\dbP_{\a^{t,l}_{s}}) +  \int_t^s \!\!\! F(r, \bar\dbP_{\a^{t,l}_r},\bar\dbP_{Y_r})dr \Big\}.
\end{equation}
where $m^*$ is the optimal argument in $\{ \bar\dbP_{Y_{t^*-}}, \bar\dbP_{Y_{t^*}} \}$, and $l^*$ the optimal argument in $\{ \bar\dbP_{\a^{t,l}_{t^*-}}, \bar\dbP_{\a^{t,l}_{t^*}} \}$. 
Clearly $t^*<T$. Indeed, if $t^* = T$, then $(T,m^*) \in \cN_{T-t}(t,m;l)$, and by \reff{c*} and \reff{utTe} we have
\beaa
&&u(T, m^*;\bar\dbP_{\a^{t,l}_T}) - \e \sum m(\dbR^d,i) +  \int_t^T F(r, \bar\dbP_{\a^{t,l}_r},\bar{\dbP}_{Y_r})dr  \ge (u-\f_\e)(t,m;l) \\
&&> \sup_{\dbP\in \cP(t, m;l)} \Big\{ \int_t^T F(r,\dbP_{\a^{t,l}_r}, \dbP_{Y_r})dr  +  H(\dbP_{\a^{t,l}_T},\dbP_{Y_T})- \sum_{\t^{\dbP}_k \geq t}\int_t^T G^k(r, \dbP_{Y_r})dr \Big\} \\
&& \ge  \int_t^T F(r, \bar\dbP_{\a^{t,l}_r}, \bar{\dbP}_{Y_r})dr + u(T, m^*;\bar\dbP_{\a^{t,l}_T}), 
\eeaa
This is the desired contradiction. Denote $ \psi^n_\e(s, \tilde m;\tilde l) :=  \f_\e(s, \tilde m) + U^{k_1,k_2} (s, \tilde m;\tilde l)$, which is obviously in $C_2^{1,2}(\overline\bQ_0)$. Here $U^{k_1,k_2}$ is in the sense of \reff{XU}, thus \reff{Uk2est} holds as well. Then, by \reff{c*},
\begin{align*}
& (u-\psi^n_\e)(t^*, m^*;l^*) =(u-\f_\e)(t^*, m^*;l^*)-U^{k_1,k_2}(t^*, m^*;l^*) \\
& =(u-\f_\e)(t^*, m^*;l^*)-\int_{t^*}^T F_{k_1}(r, \bar\dbP_{\a^{t,l}_r},\bar{\dbP}_{Y_r})dr-U^{k_1,k_2}(T, \bar\dbP_{Y_r};\bar\dbP_{\a^{t,l}_T})\\
& =(u-\f_\e)(t^*, m^*;l^*)+\big(\int_t^{t^*}-\int_t^T\big) F_{k_1}(r,\bar\dbP_{\a^{t,l}_r},\bar{\dbP}_{Y_r})dr-U^{k_1,k_2}(T,\bar\dbP_{Y_T};\bar\dbP_{\a^{t,l}_T}) \\
& \ge  (u-\f_\e)(s, \bar\dbP_{Y_s};\bar\dbP_{\a^{t,l}_s}) + \int_t^s F(r, \bar\dbP_{\a^{t,l}_r},\bar{\dbP}_{Y_r})dr + \int_t^{t^*} (F_{k_1}-F)(r, \bar\dbP_{\a^{t,l}_r},\bar{\dbP}_{Y_r})dr \\
& - \int_t^T F_{k_1}(r, \bar\dbP_{\a^{t,l}_r},\bar{\dbP}_{Y_r})dr-U^{k_1,k_2}(T,\bar\dbP_{Y_T};\bar\dbP_{\a^{t,l}_T})\\
& =  (u-\f_\e)(s, \bar\dbP_{Y_s};\bar\dbP_{\a^{t,l}_s}) - \int_s^T F_{k_1}(r, \bar\dbP_{\a^{t,l}_r},\bar{\dbP}_{Y_r})dr-U^{k_1,k_2}(T,\bar\dbP_{Y_T};\bar\dbP_{\a^{t,l}_T})+ \int_{t^*}^s (F- F_{k_1})(r, \bar\dbP_{\a^{t,l}_r},\bar{\dbP}_{Y_r})dr\\
& \ge (u-\f_\e)(s, \bar\dbP_{Y_s};\bar\dbP_{\a^{t,l}_s})-U^{k_1,k_2}(s, \bar\dbP_{Y_s};\bar\dbP_{\a^{t,l}_s})-{\e\over n} = (u-\psi^n_\e)(s, \bar\dbP_{Y_s};\bar\dbP_{\a^{t,l}_s})-{\e\over n}.
\end{align*}
Let $s=t$, by \reff{utTe} we have
\beaa
(u-\psi^n_\e)(t^*, m^*;l^*)\geq (u-\psi^n_\e)(t,m;l)-{\e\over n}>-{3\e\over n}.
\eeaa
where $(m^*;l^*)=(\bar\dbP_{Y_{t^*}};\bar\dbP_{\a^{t,l}_{t^*}})$.
Since $(u-\psi^n_\e)(T,\bar\dbP_{Y_T};\bar\dbP_{\a_T})+\e=(u-U^{k_1,k_2})(\bar\dbP_{Y_T};\bar\dbP_{\a_T}) \le (H-H_{k_1})(\bar\dbP_{Y_T},\bar\dbP_{\a_T}) \le {\e\over n}$, by compactness of $\cN_{\frac T n}(t^*, m^*;l^*)$, there exists an optimal argument $(\hat t, \hat m), \ t^*< \hat t < T$, to the above maximum. Thus $\psi^n_\e \in \ul \cA u(\hat t, \hat m;l^*)$, and therefore,
 \beaa
&&\dis  [\dbL \psi^n_\e + F + \dbQ u] (\hat t, \hat m;l^*) = \Big[\dbL \f_\e + \dbL U^{k_1,k_2} + F + \dbQ u\Big](\hat t, \hat m;l^*) \\
 &&\dis\qq =-\e + (\dbL U^{k_1,k_2} + F_{k_1} + \dbQ u)(\hat t, \hat m;l^*)  + (F-F_{k_1})(\hat t, \hat m,l^*) \\
 &&\dis\qq \le -\e + C_{L, L_{k_1}}\d_2 +  {\e\over n} \le -{\e\over 3} <0,\q ~\mbox{for large enough}~ n.
 \eeaa
This contradicts with \reff{u-contra}. The proof is completed.
\qed

\section{Two-time-scale approach}\label{sect-TTS}
\subsection{Limit system}
In this section,we assume that the Markov chain $\a^{\e}$ has a
two-time-scale structure, which is generated by $Q_{\e}=(\lambda_{pq}^{\e})$ such that
$$
Q_{\e}={1 \over \e}\tilde{Q}+\hat{Q},
$$
where $\tilde{Q}=(\tilde{\l}_{pq})$ represents the rapidly changing part and $\hat{Q}=(\hat{\l}_{pq})$ stands for the slowly moving part. Assume further that the state space of $\a^{\e}$ is given by $\cM=\cM_1\cup\cdots \cup \cM_L$, where $\cM_k =\{s_{k1},\cdots,s_{km_k}\}$ for $k= 1,\cdots,L$ and $M= m_1+\cdots+m_L$. In addition, $\tilde{Q}$ has the following block-diagonal structure
$$
\begin{pmatrix}
\tilde{Q}^1 & & & &\\
&   \cdot & & &\\
& &   \cdot & &\\
& & &  \cdot & \\
& & & &  \tilde{Q}^L
\end{pmatrix}
$$
such that $\tilde{Q}^k$ is also a generator on $\cM_k$, for $k= 1,\cdots,L$. We assume that $\tilde{Q}^k$, meaning there exists a unique  stationary distribution.

Denote $V^{\e}$ the value function of the optimal switching problem
when the time-scale parameter of the underlying Markov chain is
$\e$. We next consider the convergence of the variational inequalities
\reff{obstacle} as $\e$ goes to zero. 

\begin{lem}
\label{limit-unique}
Let $\bar{V}(\cdot,\cdot;l):=\underset{\e \rightarrow 0}{\lim}V^{\e}(\cdot,\cdot;l)$, then $\d_{l}\bar{V}(\cdot,\cdot;l,s_{kw})$ is independent of $w$, where $s_{kw} \in \cM_k =\{s_{k1},\cdots,s_{km_k}\}$.
\end{lem}
\proof
Fix $(t, m)$ and $\f\in \overline{\cA}V^{\e}(t,m;l)$, by the supersolution property,
$$
-(\dbL \f + F + \dbQ V^{\e}) (t,m;l) \ge 0.
$$
Recall $\tilde{\l}_{s_{kw}q}=0$ when $q \notin \cM_k$, multiplying both sides of the above inequality by $\e$ and letting $\e \rightarrow 0$, it follows that
$$
\sum_{r \neq w} \int_{\bS}\tilde{\lambda}_{s_{k w}, s_{k r}}(x)\left(\d_{l,m}\bar{V}(t,m,y;l, s_{k r})-\d_{l,m}\bar{V}(t,m,y;l, s_{k w})\right)m(dy)l(s_{k w}) \leq 0.
$$
That is
\begin{equation}
\label{irreduc}
\d_l\bar{V}(t,m;l, s_{k w}) \geq \sum_{r \neq w}\left(-\frac{\tilde{\lambda}_{s_{k w}, s_{k r}}}{\tilde{\lambda}_{s_{k w}, s_{k w}}}\right) \d_l\bar{V}(t,m;l,s_{k r}).
\end{equation}
By the irreducibility of $\tilde{Q}^k$ and Lemma A.39 in \cite{YZ}, the desired result is obtained. 
\qed

We remark that \eqref{irreduc} is essential for deriving limit system so that we define viscosity solutions by smoothly replacing only in $\dbL u$ instead of $\dbQ u$.
Let $\d_l\bar{u}(t,m;l, s_{k w})=\d_l\bar{u}(t,m;l, k)$, for $k= 1,\cdots,L$. We next define a limit optimal switching problem with averaged
coefficients. Let $v^k=(v^k_1,\cdots,v^k_{m_k})$  be the stationary distribution of $\tilde{Q}^k$. We define
\bea\label{limit-coeff}
\left.\ba{c}
\dis \bar{b}(k,y,m)=\sum_{w=1}^{m_k} v_w^k b(s_{k w},y,m),\q \bar{\sigma}^2(k,y,m)=\sum_{w=1}^{m_k} v_w^k \sigma^2(s_{k w},y,m),\\
\dis \bar{F}(l,m)=\sum_{k}l(k)\sum_{w=1}^{m_k} v_w^k \int_{\bS}f(s_{k w},y,m)m(dy),\q \bar{Q}=diag(v^1,\cdots,v^L)\hat{Q}diag(\1_{m_1},\cdots,\1_{m_L}). 
\ea\right. 
\eea
We note that $\bar{F}$ is no different than before, we only use it to label the limiting cases of two-time-scale Markov chain.
The corresponding limit system of variational inequalities is
\bea\label{limit-obstacle}
\dis  \min\Big\{-(\bar{\dbL} u  + \bar{F} + \bar{\dbQ} u), u - M[u]\Big\}(t,m;l) = 0, 
\eea
where 
\beaa
&&\dis \bar{\dbL} u(t,m;l)=\pa_t u(t,m;l) + \underset{k}{\sum}l(k)\underset{i \in \cN}{\sum} \int_{\dbR^d} \Big[ \bar{b}\cd \pa_x \d_{l,m} u(m,x,i;l,k) + {1\over 2} \bar{\sigma}^2: \pa_{xx}^2  \d_{l,m} u(t,m,x,i;l,k) \Big] m(dx,i)\\
&&\dis \bar\dbQ u(t,m;l):=\sum_{k}l(k)\sum_{w=1}^{m_k} v_w^k\int_{\bS}\sum_{q\neq s_{k w}}\hat\lambda_{s_{k w}q}(x)\big[\d_{l,m} u(t,m,y;l,q)-\d_{l,m} u(t,m,y;l,s_{k w})\big]m(dy)\\
&&\dis =\sum_{k}l(k)\sum_{w=1}^{m_k} v_w^k\sum_{q\neq s_{k w}}\hat\lambda_{s_{k w}q}(x)\big[\d_{l} u(t,m,y;l,q)-\d_{l} u(t,m,y;l,s_{k w})\big].
\eeaa

\begin{thm}
For $k= 1,\cdots,L$ and $w=1,\cdots,m_k$, we have $u^{\e}(t,m;l) \rightarrow \bar{u}(t,m;l)$. Moreover, $\bar{u}(t,m,k)$ is the unique viscosity solution to the system of variational inequalities \reff{limit-obstacle} for the limit optimal switching problem.
\end{thm}
\proof
Fix $(t, m)$ and $\f\in \overline{\cA}\bar{u}(t,m;l)$ with corresponding $\d$. Suppose $u^{\e}(t,m;l) \rightarrow \bar{u}(t,m;l)$, we can find $\d' < \d$ s.t. $(t_{\e},m_{\e}) \in \cN_{\d'}(t,m;l)$ is a maximizer of $\f -u^{\e}$ on $\cN_{\d'}(t,m;l)$.  We first note that
$(t_{\e},m_{\e}) \underset{\e \rightarrow 0}{\longrightarrow} (t,m)$.
Then we can find $\d'' < \d'$ s.t. $(t_{\e},m_{\e}) \in \cN_{\d'}(t,m;l)$ is a maximizer of $\f -u^{\e}$ on $\cN_{\d''}(t_{\e},m_{\e};l)$, thus $\f\in \overline{\cA}u^{\e}(t_{\e},m_{\e};l)$. Moreover,
\bea\label{obstacle-e}
\dis \sum_{k}l(k) \sum_{w=1}^{m_k}v_w^k\Big\{-\Big[\dbL u^{\e}  +  F\Big](t_{\e},m_{\e};l) - \sum_{q\neq s_{k w}}\l^{\e}_{s_{k w}q}\Big(\d_l u^{\e}(t_{\e},m_{\e};l,q)-\d_l u^{\e}(t_{\e},m_{\e};l,s_{k w})\Big) \Big\}\geq 0, 
\eea
By Lemma \ref{limit-unique}, then \reff{obstacle-e} leads to
\beaa
&&  -\Big[\dbL \bar{u}  +  F\Big](t,m;l) - \sum_{q\neq s_{k w}}\hat{\l}_{s_{k w}q}\Big(\d_l \bar{u}(t,m;l,q)-\d_l \bar{u}(t,m;l,s_{k w})\Big) \\*
&\geq& \underset{\e \rightarrow 0}{\lim} -\Big[\dbL u^{\e}  +  F\Big](t_{\e},m_{\e};l) - \sum_{q\neq s_{k w}}\l^{\e}_{s_{k w}q}\Big(\d_l u^{\e}(t_{\e},m_{\e};l,q)-\d_l u^{\e}(t_{\e},m_{\e};l,s_{k w})\Big) \geq 0. 
\eeaa
Therefore, we have
\beaa
&&  -\Big[\bar{\dbL} \bar{u}  +  \bar{F}\Big](t,m;l) - \sum_{k}l(k)\sum_{q\neq k}\bar{\l}_{kq}\Big(\d_l \bar{u}(t,m;l,q)-\d_l \bar{u}(t,m;l,k)\Big)\\
&=& \sum_{k}l(k)\sum_{w=1}^{m_k}v_w^k\Big\{-\Big[\dbL \bar{u}  +  F\Big](t,m;l) - \sum_{q\neq s_{k w}}\hat{\l}_{s_{k w}q}\Big(\d_l \bar{u}(t,m;l,q)-\d_l \bar{u}(t,m;l,s_{k w})\Big) \Big\}\geq 0. 
\eeaa
On the other hand, we also have $\Big(\bar{u}-M[\bar{u}]\Big)(t,m;l)\geq 0$ following the same analysis.

To show the subsolution property, we suppose $\Big(\bar{u}-M[\bar{u}]\Big)(t,m;l)> 0$. $\forall \f\in \underline{\cA}\bar{u}(t,m;l)$ with corresponding $\d$. We can find $\d' < \d$ s.t. $(t_{\e},m_{\e}) \in \cN_{\d'}(t,m;l)$ is a minimizer of $\f -u^{\e}$ on $\cN_{\d'}(t,m;l)$.  We first note that
$(t_{\e},m_{\e}) \underset{\e \rightarrow 0}{\longrightarrow} (t,m)$.
Then we can find $\d'' < \d'$ s.t. $(t_{\e},m_{\e}) \in \cN_{\d'}(t,m;l)$ is a minimizer of $\f -u^{\e}$ on $\cN_{\d''}(t_{\e},m_{\e};l)$, thus $\f\in \underline{\cA}u^{\e}(t_{\e},m_{\e};l)$. For $\e$ small enough, we have $\Big(u^{\e}-M[u^{\e}]\Big)(t_{\e},m_{\e};l)> 0$. It follows that
\beaa
\dis \sum_{k}l(k) \sum_{w=1}^{m_k}v_w^k\Big\{-\Big[\dbL u^{\e}  +  F\Big](t_{\e},m_{\e};l) - \sum_{q\neq s_{k w}}\l^{\e}_{s_{k w}q}\Big(\d_l u^{\e}(t_{\e},m_{\e};l,q)-\d_l u^{\e}(t_{\e},m_{\e};l,s_{k w})\Big)\Big\}\leq 0, 
\eeaa
then we can repeat the same argument in supersolution part and obtain the required result.
\qed
\subsection{Examples}
In this subsection we construct examples where the obstacle problem indeed has a viscosity solution. For simplicity, we consider a two-mode case, i.e. $\cN=\{0,1\}$. We first start with a two-state Markov chain $\a$, i.e. $\cM=\{1,2\}$. The generator $Q$ is,
$$
\begin{pmatrix}
-\m_1 & \m_1\\
\m_2  & -\m_2
\end{pmatrix}
$$
 Let $\psi\in C( \dbR)$, $a^i(q) \in C([0, T])$, $i\in\cN, q\in\cM$ be positive functions, and set
\bea
 u(t, m;l) := (T-t)\underset{q\in\cM}{\sum}l(q)\big(\big[v^1(m)\wedge a^1(q)\big]+\big[v^0(m)\wedge a^0(q)\big]\big),\q\mbox{where}\q v^i(m) := \int_\dbR \psi(x) m(dx,i).
\eea

\begin{prop}
\label{prop-eg}
Under the above setting, $u \in C([0, T]\times \cP_2(\bf S\times \cM))$, and $u$ is the viscosity solution to the obstacle problem \reff{obstacle} with
\bea
\label{data}
&&\dis F(t,m;l) =\underset{q\in\cM}{\sum}l(q)\big\{|v^1(m)-a^1(q)|+|v^0(m)-a^0(q)|+a^1(q)+a^0(q)\big)\big\},\\
&&\dis \m_1=\cL^{i,1},\q \m_2=\cL^{i,2}, \q g_{10}(t,x) =  g_{01}(t,x) = (T-t)\psi(x), \q H:=0. \nonumber
\eea
\end{prop} 

\proof  First, it's easy to check $u \in C([0, T]\times \cP_2(\bf S))$ and $g_{10}(t,x)+g_{01}(t,x)>0$.
We now show that $u$ satisfies \reff{obstacle}. Clearly, $u(T, \cdot) = 0 = H$, and
\beaa
&&\dis \dbL u(t,m;l) = \underset{q\in\cM}{\sum}l(q)\big\{-\big[v^0(m)\wedge a^0(q)\big]-\big[v^1(m)\wedge a^1(q)\big]\\
&&\dis +(T-t)\int_{\dbR}\big( \cL^{0,q}_x\d_m \big[v^0(m)\wedge a^0(q)\big] m(dx,0)+\cL^{1,q}_x\d_m \big[v^1(m)\wedge a^1(q)\big] m(dx,1)\big)\big\},\\
&&\dis \dbQ u(t,m;l)=(T-t)(\m_1 l(1)-\m_2 l(2))\Big\{ \int_{\dbR}\big(\d_m \big[v^0(m)\wedge a^0(2)-v^0(m)\wedge a^0(1)\big] m(dx,0)\\
&&\dis +\d_m \big[v^1(m)\wedge a^1(2)-v^1(m)\wedge a^1(1)\big] m(dx,1)\big)\Big\}.
\eeaa
Then by $\cL^{i,1}=\m_1,\cL^{i,2}=\m_2$, we have
\bea
-(\dbL u  + F + \dbQ u)(t,m;l)=\underset{q\in\cM}{\sum}l(q)\big\{ [ v^1(m)-a^1(q)]^+ +[ v^0(m)-a^0(q)]^+\big\}\geq 0.
\eea
Thus, the beginning and ending of a switching should be 
\bea\label{region}
\{v^0(m)>a^0(q)\}\cup\{v^1(m)>a^1(q)\}\longrightarrow \{v^0(m')\leq a^0(q)\}\cap \{v^1(m')\leq a^1(q)\}.
\eea 
Finally, we need to verify $g_{10}(t,x), g_{01}(t,x)$ satisfy $\big(u - M[u]\big)(t,m,q)\geq 0$. 
Note that $\int_{\dbR}\psi(x)p_{11}m(dx,1)=v^1(m')$ for a switching from $1$ to $0$ and $\int_{\dbR}\psi(x)p_{00}m(dx,0)=v^0(m')$ for $0$ to $1$, set $v^1(m'):=x$ and we have
\beaa
\frac{u - M[u]}{T-t}(t,m;l)=\begin{cases}
  \underset{q\in\cM}{\sum}l(q)\Big\{\big[v^1(m)\wedge a^1(q)\big]-[v^0(m)-a^0(q)]^+ -\underset{m' \prec m} {\sup}\big[x\wedge a^1(q)-[v^0(m')-a^0(q)]^+\big]\Big\},\q  x<v^1(m)\\
  \underset{q\in\cM}{\sum}l(q)\Big\{\big[v^0(m)\wedge a^0(q)\big]-[v^1(m)-a^1(q)]^+ -\underset{m' \prec m} {\sup}\big[v^0(m')\wedge a^0(q)-[x-a^1(q)]^+\big]\Big\},\q  x>v^1(m)\\
\end{cases}
\eeaa
Thus, $\big(u - M[u]\big)(t,m,q)\geq 0$, and "=" holds if and only if $\{v^1(m)>v^1(m')\geq a^1(q)\}\cap\{v^0(m)<v^0(m')\leq a^0(q)\}$ and $\{v^0(m)>v^0(m')\geq a^0(q)\}\cap\{v^1(m)<v^1(m')\leq a^1(q)\}$. By \reff{region}, the switching is required to be
$$
\{v^0(m)>a^0(q)\}\longrightarrow \{v^0(m')= a^0(q)\}\cap\{v^1(m')\leq a^1(q)\},\q \{v^1(m)>a^1(q)\}\longrightarrow \{v^1(m')= a^1(1)\}\cap\{v^0(m')\leq a^0(q)\}.
$$
The proof is completed.
\qed

If we view this example as stocks trading: let the $i=1$  represent "long sale" and $i=0$  "short sale". Both long and short can be gained, one needs to choose the optimal position in shares $m(1)$ and $m(0)$, i.e. long-short optimal capital allocation, according to the value function $u(t,m)$. We emphasize that the strategy here allows a distribution of $m(1)$ and $m(0)$, which is exactly the significant difference between mean field optimal switching and the standard case. The coefficients in \reff{data} imply that $F$, i.e holding gains are influenced by the long-short capital allocation and market status.

Let $q=2$ represent "bull market" and $q=1$ "bear market", then the profit-taking price $a^1(2)>a^1(1), a^0(2)<a^0(1)$. According to \reff{region}, $\{v^1(m)>a^1(q)\}$ is the region one needs to switch short  and $\{v^0(m)>a^0(q)\}$ switch long. As shown in  Table \ref{tab1}, the "long" region is larger and "short" region is smaller when $q=2$. 
In Table \ref{tab1}, the black letters represent switching strategy while yellow letters represent preserving original strategy. For example, one needs to keep short in "bear market" ($q=1$) while switch long in "bull market" ($q=2$) in region \uppercase\expandafter{\romannumeral3}. In real life, we know  whether we are currently in a bull or bear market,  but our framework allows for being in a bull or bear market with some probability.

\begin{figure}[!htb]
  \centering
  \includegraphics[width=0.7\linewidth]{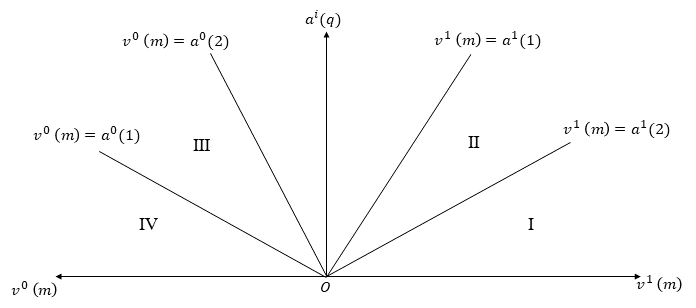}
  \caption{Switching region for $q=1, 2$}
  \label{fig1}
\end{figure}

\begin{table}[!ht]
  \centering
  \caption{Switching strategy for $q=1, 2$}
  \label{tab1}
  \begin{tabular}{cccccc}
  \toprule
  \multicolumn{2}{c}{region}  &  \uppercase\expandafter{\romannumeral1}   & \uppercase\expandafter{\romannumeral2}    & \uppercase\expandafter{\romannumeral3}   & \uppercase\expandafter{\romannumeral4} \\
  \midrule
  strategy & q=1 & short & short & \textcolor{yellow}{short} & long \\
           & q=2 & short &  \textcolor{yellow}{long} & long  & long \\
  \bottomrule
  \end{tabular}
  \end{table}

In the rest of this subsection, we construct an example to explain how two-time-scale approach works. We further consider a four-state case with Markov chain $\a(t) = (\a_1(t), \a_2(t))$, where $\a_1(t) \in \{1, 2\}$ represents the primary market trend  and $\a_2(t) \in \{1, 2\}$ represents the secondary market movement indicator. The generator $Q^{\e}$ is,
$$
\frac{1}{\e}
\begin{pmatrix}
-\l_1 & \l_1 & 0 & 0\\
\l_2  & -\l_2 & 0 & 0\\
0 & 0 & -\l_1 & \l_1\\
0 & 0 & \l_2  & -\l_2
\end{pmatrix}
+
\begin{pmatrix}
-\m_1 & 0 & \m_1 & 0\\
0 &  -\m_1 & 0 & \m_1\\
\m_2 & 0 & -\m_2 & 0\\
0 & \m_2 & 0 & -\m_2
\end{pmatrix}
$$ 

We also denote as $\a^{\e}(t)\in \{1, 2, 3, 4\}$. In the below, the state $q \in \{1, 2, 3, 4\}$ refers to $(q_1,q_2)\in \{(1,1),(1,2),(2,1),(2,2)\}$ respectively. Let  $a^i(q_1,q_2) \in C([0, T])$ be positive functions, and set
\beaa
 u^{\e}\big(t, m;l\big) &:=& (T-t)\underset{(q_1,q_2)}{\sum}l(q_1,q_2)\Big\{\big[v^1(m)\wedge a^1(q_1,q_2)\big]+\big[v^0(m)\wedge a^0(q_1,q_2)\big]\Big\}.
\eeaa

\begin{prop}
\label{prop-eg2}
Under the above setting, $u^{\e} \in C([0, T]\times \cP_2(\bf S))$, and $u^{\e} $ is the viscosity solution to the obstacle problem \reff{obstacle} with
\beaa
&&\dis F(t,m;l) =\underset{(q_1,q_2)}{\sum}l(q_1,q_2)\big\{|v^1(m)-a^1(q_1,q_2)|+|v^0(m)-a^0(q_1,q_2)|+a^1(q_1,q_2)+a^0(q_1,q_2)\big)\big\},\\
&&\dis  g_{10}(t,x)=g_{01}(t,x) = (T-t)\psi(x),\q \m_1=\cL^{i,(1,q_2)},\q \m_2=\cL^{i,(2,q_2)}, \q H:=0.
\eeaa
$\bar{u}(t,,m,k)= (T-t)\underset{k=1,2}{\sum}l(k)\Big\{\big[v^1(m)\wedge a^1(k)\big]+\big[v^0(m)\wedge a^0(k)\big]\Big\}$ is the viscosity solution of \reff{obstacle} with 
\beaa
 \bar{F}(t,m,k) = \underset{k=1,2}{\sum}l(k)\big\{|v^1(m)-a^1(k)|+|v^0(m)-a^0(k)|+a^1(k)+a^0(k)\big)\big\}.
\eeaa
\end{prop} 

\proof  First, it's easy to check $u^{\e} \in C([0, T]\times \cP_2(\bf S \times \cM))$ and $g_{10}(t,x)+g_{01}(t,x)>0$.
We now show $u^{\e}$ satisfies \reff{obstacle}. Clearly, $u^{\e}(T, \cdot) = H$, and let $q \in \{1, 2, 3, 4\}$ refers to $(q_1,q_2)\in \{(1,1),(1,2),(2,1),(2,2)\}$ respectively .
\beaa
&&\dis \dbL u^{\e}(t,m,q) = \underset{q\in\cM}{\sum}l(q)\Big\{-\big[v^1(m)\wedge a^1(q)\big]-\big[v^0(m)\wedge a^0(q)\big]\\
&&\dis +(T-t)\int_{\dbR} \big(\cL^{1,q}_x\d_{l,m} u^\e(
t,m,x,1;l,q) m(dx,1)  +\cL^{0,q}_x\d_{l,m} u^\e(
t,m,x,0;l,q) \big) m(dx,0) \Big\},\\
&&\dis \dbQ u^\e(t,m;l)=(T-t)\underset{p\in\cM}{\sum}l(p)\int_{\bS}\underset{q\neq p}{\sum}\l_{pq}^{\e}[\d_{l,m}u^\e(
t,m,y;l,q)-\d_{l,m}u^\e(
t,m,y;l,p)]m(dy)\\
&&\dis=(T-t)\Big\{(\m_1 l(1)-\m_2 l(3))\int_{\bS}[\d_{l,m}u^\e(
t,m,y;l,3)-\d_{l,m}u^\e(t,m,y;l,1)]m(dy)\\
&&\dis + (\m_1 l(2)-\m_2 l(4))\int_{\bS}[\d_{l,m}u^\e(
t,m,y;l,4)-\d_{l,m}u^\e(t,m,y;l,2)]m(dy)\\
&&\dis + \frac{1}{\e}(\l_1 l(1)-\l_2 l(2))\int_{\bS}[\d_{l,m}u^\e(t,m,y;l,2)-\d_{l,m}u^\e(t,m,y;l,1)]m(dy)\\*
&&\dis + \frac{1}{\e}(\l_1 l(3)-\l_2 l(4))\int_{\bS}[\d_{l,m}u^\e(t,m,y;l,4)-\d_{l,m}u^\e(t,m,y;l,3)]m(dy)\Big\}.
\eeaa
Note that $\m_1=\cL^{i,(1,q_2)},\q \m_2=\cL^{i,(2,q_2)}$,  then we have 
\small{
\bea
\label{diff-con}
&&\dis -(\dbL u^{\e}  + F + \dbQ u^{\e})(t,m;l)=\underset{q\in\cM}{\sum}l(q)\big\{ [ v^1(m)-a^1(q)]^+ +[ v^0(m)-a^0(q)]^+\big\} - \frac{1}{\e}(T-t)\Big\{ \\
&&\dis (\l_1 l(1)-\l_2 l(2))[\d_{l}u^\e(m,y;l,2)-\d_{l}u^\e(m,y;l,1)]+(\l_1 l(3)-\l_2 l(4))[\d_{l}u^\e(m,y;l,4)-\d_{l}u^\e(m,y;l,3)]\Big\} \nonumber.
\eea
}
Note the stationary distributions require $\frac{l(1)}{l(2)}=\frac{l(3)}{l(4)}=\frac{\l_2}{\l_1}$ as $\e\rightarrow 0$,
then \reff{diff-con} equals to
\bea\label{diff-reg}
-(\dbL u^{\e}  + F + \dbQ u^{\e})(t,m;L)=\underset{q\in\cM}{\sum}l(q)\big\{ [ v^1(m)-a^1(q)]^+ +[ v^0(m)-a^0(q)]^+\Big\}\geq 0,\q \e\rightarrow 0.
\eea
By the continuity, we obtain the convergence of diffusion condition. It should be noted that the optimal strategy depends both on $a(q)$ and $\e$ for given $\e$, while the strategy for limit system only depends on $a(q)$.

Finally, we can verify $g_{10}(t,x), g_{01}(t,x)$ satisfy $\big(u^{\e} - M[u^{\e}]\big)(t,m,q)\geq 0$ same with proposition
 \ref{prop-eg}. Moreover, the switching reads
\beaa
\{v^0(m)>a^0(q_1,q_2)\}\longrightarrow \{v^0(m')= a^0(q_1,q_2)\}\cap\{v^1(m')\leq a^1(q_1,q_2)\},\\
\{v^1(m)>a^1(q_1,q_2)\}\longrightarrow \{v^1(m')= a^1(q_1,q_2)\}\cap\{v^0(m')\leq a^0(q_1,q_2)\}.
\eeaa
Thus, we obtain the convergence of switching condition. Note that $a^i(q_1,q_2)\rightarrow a^i(k)$, $k=q_1$, the convergence of $u^\e$ can be easily verified. The proof is completed.
\qed

\begin{figure}[!htb]
  \centering
  \includegraphics[width=0.7\linewidth]{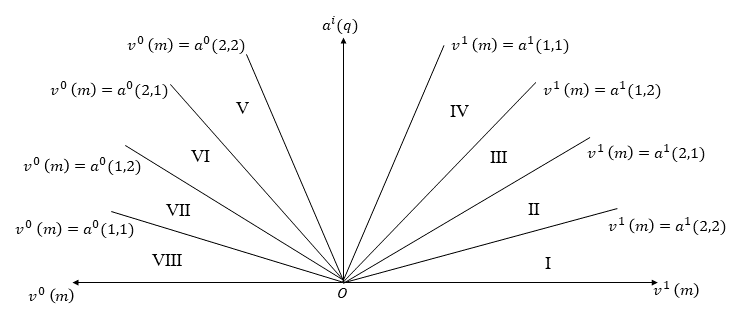}
  \caption{Switching region for $q=1, 2, 3, 4$}
  \label{fig2}
\end{figure}

In this example, let $q=2$ represent "bull market" and $q=1$ "bear market", then the profit-taking price $a^1(2,2)\geq a^1(2,1)\geq a^1(1,2)\geq a^1(1,1)$ and $a^0(2,2)\leq a^0(2,1)\leq a^0(1,2)\leq a^0(1,1)$.  As shown in Figure \ref{fig2} and Table \ref{tab2}, the "long" region becomes larger and "short" region smaller with the increasing of $q$.  In Table \ref{tab2}, the black letters represent switching strategy while yellow letters represent preserving original strategy. For limit system, there must be $a^i(q_1,q_2)=a^i(k)$, $k=q_1$ as $\e \to 0$. Then we can treat it as a two Markov state case with $k=1, 2$, i.e. the problem in Proposition \ref{prop-eg}. It's not difficult to observe that the strategy shown in Figure \ref{fig2} and Table \ref{tab2} transforms into Figure \ref{fig1} and Table \ref{tab1}. By replacing with limit coefficients in (\ref{limit-coeff}), the computational procedure does get simplified.

\begin{table}[!ht]
  \centering
  \caption{Switching strategy for $q=1, 2, 3, 4$}
  \label{tab2}
  \scalebox{0.7}{
  \begin{tabular}{cccccccccc}
  \toprule
  \multicolumn{2}{c}{region}   & \uppercase\expandafter{\romannumeral1}   & \uppercase\expandafter{\romannumeral2}    & \uppercase\expandafter{\romannumeral3} & \uppercase\expandafter{\romannumeral4}   & \uppercase\expandafter{\romannumeral5}    & \uppercase\expandafter{\romannumeral6} & \uppercase\expandafter{\romannumeral7}   & \uppercase\expandafter{\romannumeral8}     \\
  \midrule
  strategy & q=1 & short & short & short & short & \textcolor{yellow}{short} & \textcolor{yellow}{short} & \textcolor{yellow}{short} & long  \\
           & q=2 & short & short & short & \textcolor{yellow}{long} & \textcolor{yellow}{short} & \textcolor{yellow}{short} & long & long  \\
           & q=3 & short & short & \textcolor{yellow}{long} & \textcolor{yellow}{long} & \textcolor{yellow}{short} & long & long & long  \\
           & q=4 & short & \textcolor{yellow}{long} & \textcolor{yellow}{long} & \textcolor{yellow}{long} & long & long & long & long \\
  \bottomrule
  \end{tabular}
  }
  \end{table}

\section{Concluding remarks}\label{sect-conc}

In this paper, we formulate a general mean field optimal switching problem under a regime-switching model and prove that the value function is indeed the unique viscosity solution of the corresponding variational inequality on the Wasserstein space. Specifically, we construct the DPP and select the compact neighborhood based on the conditional distribution of Markov chain. To overcome the computational difficulty brought by the appearance of Markov chain, we employ the so-called two-time-scale approach  and obtain a convergence result for limit system. Moreover, we provide an innovative and interpretive example,  which is further used in the massive stock trading problem.

However, it is a pity that the example is just for illustration purpose of the theory, for the computation of the obstacle problem on the Wasserstein space remains a challenge task. In future, we shall focus on providing an explicit calculation algorithm by exploiting the finite dimensional approximation. We may estimate the value function with a neural network by choosing appropriate loss function.

\appendix

\section{Proof of Proposition \ref{Itothm}}
\label{appendixA}

 Denote $\D Y_s:= Y_s-Y_{s-}$ and $Y^D_t := \sum_{0<s\le t} \D Y_s$. By \reff{summability},
\bea
\label{EY*}
\dbE\big[|Y^*_T|^2 + \|Y^D\|_T^2\big]<\infty,\q\mbox{where}\q Y^*_T := \sup_{0\le s\le T} |Y_s|,\q  \|Y^D\|_t := \sum_{0<s\le t} |\D Y_s|.
\eea
 For $n\ge 1$, set $\D t := {T\over n}$, $t_i:= i \D t$, $i=0,\cds, n$., then 
\begin{equation}\label{Du}
\begin{split}
&u(t_{i+1}, m_{t_{i+1}};l_{t_{i+1}})-u(t_i, m_{t_i};l_{t_i})\\
&=u(t_{i+1}, m_{t_{i+1}};l_{t_{i+1}})-u(t_i, m_{t_{i+1}};l_{t_{i+1}})+u(t_i, m_{t_{i+1}};l_{t_{i+1}})-u(t_i, m_{t_i};l_{t_i})\\
& =\int_{t_i}^{t_{i+1}} \pa_t u(s, m_{t_{i+1}};l_{t_{i+1}}) ds + \int_0^1 \dbE[\xi^\l_{t_{i}}] d\l .
\end{split}
\end{equation}
Here $\xi^\l_{t_{i}}:= \d_{l,m} u(t_i, m_{t_i}^\l,Y_{t_{i+1}};l^{\l}_{t_i},\a_{t_{i+1}})-\d_{l,m} u(t_i, m_{t_i}^\l,Y_{t_i};l^{\l}_{t_i},\a_{t_i})$,  $\psi^\l_{t_i} := \l \psi_{t_i} + [1-\l]\psi_{t_{i+1}}, \psi=l,m$.
By the standard It\^{o}'s formula:
\beaa
&\dis  \xi^\l_{t_i} = \d_{l,m} u(t_i, m_{t_i}^\l,Y_{t_{i+1}};l^{\l}_{t_i},\a_{t_{i+1}})-\d_{l,m} u(t_i, m_{t_i}^\l,Y_{t_{i+1}};l^{\l}_{t_i},\a_{t_i})+\d_{l,m} u(t_i, m_{t_i}^\l,Y_{t_{i+1}};l^{\l}_{t_i},\a_{t_i})-\d_{l,m} u(t_i, m_{t_i}^\l,Y_{t_i};l^{\l}_{t_i},\a_{t_i})\\
&\dis =\int_{t_i}^{t_{i+1}}  \big[\G^{1,\l}_s+\G^{2,\l}_s\cd dY_s^c  +  \frac{1}{2} \G^{3,\l}_s:d\langle Y^c \rangle_s\big] ds + \int_{(t_i, t_{i+1}]} \G^{4,\l}_s dY^D_s,\\
&\dis \mbox{where} \q \G^{1,\l}_s:=Q \d_{l,m} u(t_i, m^{\l}_{t_i},Y_{t_{i+1}};l_{t_i}^\l,\a_s) ,\q \G^{2,\l}_s :=  \pa_y \d_{l,m} u(t_i, m_{t_i}^\l,Y_s;l^{\l}_{t_i},\a_{t_i}),\\ 
&\dis \G^{3,\l}_s :=  \pa_{yy}^2 \d_{l,m} u(t_i, m_{t_i}^\l,Y_s;l^{\l}_{t_i},\a_{t_i}),\q \G^{4,\l}_s:= \int_0^1 \pa_y \d_{l,m} u\big(t_i, m_{t_i}^\l,\th Y_s + [1-\th] Y_{s-};l^{\l}_{t_i},\a_{t_i}\big) d\th.
\eeaa
Note that $m_{t_{i+1}}, m^\l_{t_i}\in \Xi_\bm$, by the growth conditions, we have 
\bea
\label{Gammabound}
 |\G^1_s|\le C,\q |\G^{2,\l}_s| \le C[1+|Y_s|],\q  |\G^{3,\l}_s|\le C,\q  |\G^{4,\l}_s|  \le C[1+|Y_s|+|Y_{s-}|].
\eea
Then
\beaa
\dbE\Big[ \Big(\int_{t_i}^{t_{i+1}}\G^{2,\l}_s(\G^{2,\l}_s)^\top  : d\la M^c\ra_s\Big)^{1\over 2}\Big] \le C\dbE\Big[ \Big([1+|Y^*_T|^2] \la M^c\ra_T\Big)^{1\over 2}\Big] \le C\dbE\Big[ 1+|Y^*_T|^2 +\la M^c\ra_T\Big]<\infty.
\eeaa
This implies $\int_0^1 \dbE\big[\int_{t_i}^{t_{i+1}} \G^{2,\l}_s \cd dM^c_s\big]d\l=0$. Sum for \reff{Du},
\beaa
&&\dis u(T, m_{T-};l_{T-}) = u(0, m_{0-};l_{0-}) +  \int_0^T \pa_t u(s, m_s;l_s) ds \\
&&\dis + \int_0^1\dbE\Big[\int_0^T \big[\G^1_s+\G^{2,\l}_s \cd d A^c_s + \G^{3,\l}_s : d\la Y^c\ra_s\big]ds + \int_{(0, T]} \G^{4,\l}_s dY^D_s\Big] d\l.
\eeaa
Fix $\l$, $s$, and send $n\to \infty$, by the regularity of $u$ we have:  denoting $m^\l_s:= \l m_{s-} + [1-\l] m_s$,
\beaa
&\G^1_s \to Q \d_{l,m} u(s, m_s,Y_s;l_s,\a_s),\q \G^{2,\l}_s\to   \pa_y \d_{l,m} u(s, m^\l_s,Y_s;l_s,\a_s),\q \G^{3,\l}_s \to  \pa_{yy}^2 \d_{l,m} u(s, m^\l_s,Y_s;l_s,\a_s),\\
& \G^{4,\l}_s\to \int_0^1 \pa_y \d_m u\big(s, m_s^\l, \th Y_s + [1-\th] Y_{s-};l_s,\a_s\big) d\th,\q a.s.
 \eeaa
By \reff{summability}, \reff{EY*}, and \reff{Gammabound},  we may apply the dominated convergence theorem to obtain
\beaa
&& u(T, m_{T-};l_{T-}) = u(0, m_{0-};l_{0-}) + \int_0^T \pa_t u(s, m_s;l_s) ds \\
&& + \int_0^T \sum_{p\in\cM}l_s(p)\sum_{q\neq p}\int_{\bS}\lambda_{pq}(x)\big[\d_{l,m} u(s, m_s,y;l_s,q)-\d_{l,m} u(s, m_s,y;l_s,p)\big]m_s(dy) ds \\
&& + \int_0^1\dbE\Big[\int_0^T \big[ \pa_y \d_{l,m} u(s, m^\l_s,Y_s;l_s,\a_s) \cd d A^c_s +  \pa_{yy}^2 \d_{l,m} u(s, m^\l_s,Y_s;l_s,\a_s) : d\la Y^c\ra_s\big]\\
&& + \int_{(0, T]} \int_0^1 \pa_y \d_{l,m} u\big(s, m_s^\l, \th Y_s + [1-\th] Y_{s-};l_s,\a_s\big) d\th dY^D_s\Big] d\l.
\eeaa
Since $A^c$ and $\la Y^c\ra$ are continuous, and $ m_s$ has at most countably many jumps, then
\bea
\label{Ito1}
&& u(T, m_{T-};l_{T-}) = u(0, m_{0-};l_{0-})  +  \int_0^T \Big[\pa_t u(s, m_s;l_s) + \dbQ u(s, m_s;l_s)\Big]  ds  \nonumber\\*
&&\q +\dbE \int_0^T\!\!\! \big[ \pa_y \d_{l,m} u(s, m_s,Y_s;l_s,\a_s) \cd d A^c_s+  \pa_{yy}^2 \d_{l,m} u(s, m_s,Y_s;l_s,\a_s) : d\la Y^c\ra_s\big]+  \cJ_D , \\*
&&\mbox{where}\q \cJ_D:=  \dbE\Big[ \int_0^1 \int_{(0, T]} \int_0^1 \pa_y \d_{l,m} u\big(s, m_s^\l, \th Y_s + [1-\th] Y_{s-};l_s,\a_s\big) d\th dY^D_s d\l\Big] .\nonumber
\eea
It remains to compute $\cJ_D$. First, by Fubini's theorem,
\bea
\label{JD0}
&\dis\cJ_D =  \dbE\Big[ \sum_{s\in (0, T]}\int_0^1  \int_0^1 \pa_y \d_{l,m} u\big(s, m_s^\l, \th Y_s + [1-\th] Y_{s-};l_s,\a_s\big) d\th  d\l \D Y_s\Big] =\dbE\Big[ \sum_{s\in (0, T]} \D \d_{l,m} u_s\Big]\nonumber\\
&\dis\mbox{where}\q \D \d_{l,m} u_s:= \int_0^1 \big[ \d_{l,m} u(s, m_s^\l, Y_s;l_s,\a_s ) - \d_{l,m} u(s, m_s^\l, Y_{s-};l_s,\a_s )\big]  d\l.
\eea
Note that $(0, T] = J_{(0, T]}(\bm) \cup J^c_{(0, T]}(\bm)$. Since $J_{(0, T]}(\bm)$ is countable, then
\bea
\label{JD1}
\dbE\Big[ \!\!\! \sum_{s\in J_{(0, T]}(\bm) } \!\!\! \D \d_{l,m} u_s\Big]=  \!\!\!\sum_{s\in J_{(0, T]}(\bm) } \!\!\!\dbE\big[ \D \d_{l,m} u_s\big] =  \!\!\! \sum_{s\in J_{(0, T]}(\bm) } \!\!\! \big[ u(s, m_s;l_s) - u(s, m_{s-};l_s)\big].
\eea
For $s\in J^c_{(0, T]}(\bm)$, we have $m^\l_s = m_s$, $0\le \l\le 1$. Then
\bea
\label{JD2}
\dbE\Big[ \!\!\! \sum_{s\in J^c_{(0, T]}(\bm) } \!\!\! \D \d_{l,m} u_s\Big]= \dbE \Big[ \!\!\! \sum_{s\in J^c_{(0, T]}(\bm) } \!\!\! \big[ \d_{m} u(s, m_s,Y_s;l_s) - \d_{m} u(s, m_s, Y_{s-};l_s)\big]\Big].
\eea
We emphasize that since $J^c_{(0, T]}(\bm)$ is uncountable, unlike in \reff{JD1} we cannot switch the order of $\dbE$ and $\sum_{s\in J^c_{(0, T]}(\bm) }$ at above. Plug \reff{JD1}, \reff{JD2} into \reff{JD0}, then plug \reff{JD0} into \reff{Ito1}, we complete the proof. \qed

\section{Proof of Lemma \ref{thm-reg}} 
\label{appendixB}

Let $\rho_0$ denote the modulus of continuity of $f, h$ under $\cW_1$. We proceed in two steps.

{\it Step 1.} Fix $t\in [0, T]$ and $m, \tilde m\in \cP_2(\bS)$. For any $\dbP\in \cP(t,m;l)$, by possibly enlarging the space, there exists $(\tilde X_t, \tilde I_{t-})$ on the space  $(\O, \cF, \dbP)$ such that 
\beaa
\dbP_{(\tilde X_t, \tilde I_{t-})} = \tilde m,~ \dbE^\dbP\Big[| \tilde X_t - X_t|+ |\tilde I_{t-} - I_{t-}|\Big] = \cW_1(m, \tilde m).
\eeaa
 Consider the following SDE on the space $(\O, \cF, \dbP)$:  for $\tilde Y := (\tilde X, \tilde I)$,
 \beaa
\label{tildeX1}
\tilde X_s = \tilde X_t + \int_t^s b(r,\a^\dbP_r, \tilde X_r, \tilde I_r, \dbP_{\tilde Y_r}) dr + \int_t^s \sigma(r, \a^\dbP_r, \tilde X_r, \tilde I_r, \dbP_{\tilde{Y}_r}) dW_r^\dbP, \q \dbP\mbox{-a.s.}\\
\tilde I_r= \underset{\t_n^{\dbP} \geq t}{\sum}\xi_{n-1}\1_t^n ~\mbox{with} ~ \xi_{\hat n -1}=\tilde I_{t-}, \q I_r= \underset{\t_n^{\dbP} \geq t}{\sum}\xi_{n-1}\1_t^n ~\mbox{with} ~ \xi_{\hat n -1}=I_{t-}.
\eeaa
Here $\hat n := \inf\{n, \t_n^{\dbP} \geq t\}$. Denote  $\D Y := \tilde Y-Y$, then 
\bea
\label{DI}
\sup_{t\le r\le T} |\D I_r| \le |\D I_{t-}|,\q\mbox{and thus}\q \dbE^\dbP\Big[\sup_{t\le r\le T} |\D I_r|\Big] \le  \cW_1(m, \tilde m).
\eea
Moreover,  for $\f=b, \si$, by the desired Lipschitz continuity under $\cW_1$, we have
\beaa
\Big|\f(r, \a^\dbP_r, \tilde X_r, \tilde I_r,\dbP_{\tilde Y_r)}) - \f(r, \a^\dbP_r, X_r, I_r, \dbP_{Y_r}) \Big|  \le C\Big[|\D X_r| + \cW_1(\dbP_{\tilde Y_r}, \dbP_{Y_r}) + |\D I_{t-}|\Big].
\eeaa
By standard estimates, one can show that
\bea
&&\dis \dbE^\dbP_t\Big[\sup_{t\le s\le T} |X_s|^2\Big] \le C[1+|X_t|^2];\nonumber\\
&&\dis \dbE^\dbP_t\Big[|\D X_s|^2\Big] \le C\int_s^T \cW_1^2(\dbP_{\tilde Y_r}, \dbP_{Y_r})dr  + C|\D X_t|^2+C\sup_{t\le s\le T}\dbE^\dbP_t |\D I_{t-}|^2;\nonumber\\
\label{DX}
&&\dbE^\dbP_t\big[ |\D X_s|\big] \le C\Big(\int_s^T \cW_1^2(\dbP_{\tilde Y_r}, \dbP_{Y_r})dr\Big)^{1\over 2} +C\big(|\D X_t|+ |\D I_{t-}|\big).
\eea
This implies that
\beaa
 \cW_1^2(\dbP_{\tilde Y_s}, \dbP_{Y_s}) \le C\int_s^T \cW_1^2(\dbP_{\tilde Y_r}, \dbP_{Y_r})dr + C\Big(\dbE^\dbP\Big[|\D X_t|+ |\D I_{t-}| \Big]\Big)^2.
\eeaa
By Grownwall inequality we have, for any $R>0$,
 \bea
 \label{dR}
 \sup_{t\le s\le T} \cW_1(\dbP_{\tilde Y_s}, \dbP_{Y_s}) \le C \cW_1(m, \tilde m) .
 \eea
Notice that $\tilde \dbP := \dbP \circ (\tilde X, I)^{-1} \in \cP(t, \tilde m;l)$. Then 
\beaa
&&\dis V(t, m;l) - V(t, \tilde m;l)\\
&&\dis \le \dbE^\dbP \int_t^T \big[f(r, \a^\dbP_r, X_r, I_r, \dbP_{Y_r})  - f(r, \a^\dbP_r, \tilde X_r, \tilde I_r, \dbP_{\tilde Y_r}) \big]dr + \big[H(\dbP_{\a^{t,l}_T},\dbP_{Y_T})  - H(\dbP_{\a^{t,l}_T},\dbP_{\tilde Y_T}) \big]\\
&&\dis -\underset{\t_n^{\dbP} \geq t}{\sum}\big[\int_t^T G^n(r, \dbP_{Y_r})dr- \int_t^T G^n(r, \dbP_{\tilde Y_r})dr\big]\\
&&\dis \le  \rho_0\Big( \cW_1(\dbP_{\tilde Y_T}, \dbP_{Y_T})\Big) + \int_t^T \dbE^\dbP\Big[\rho_0\big(|\D Y_r|)+ \rho_0\big(\cW_1(\dbP_{\tilde Y_r}, \dbP_{Y_r})\big) \Big]dr \le C_m \cW_1(m,\tilde m) + \int_t^T \dbE^\dbP\big[\rho_0(|\D X_r|)] dr.
\eeaa
Fix $m$ and send $\tilde m\to m$ under $\cW_1$, we see that 
 $\underset{\tilde m\to m}{\lim} V(t, \tilde m;l) = V(t,m;l)$.

{\it Step 2.} Let $t<\tilde t$ and $(m,l)\in \cP_2(\bS\times\cM)$.   By DPP we have
\bea
\label{A-DPP}
\left.\ba{c}
\dis V(t, m;l) = \sup_{\dbP\in \cP(t, m;l)} \Big\{ \int_t^{\tilde t} \!\! F(r,  \dbP_{\a^{t,l}_r}, \dbP_{Y_r}) dr -\underset{\t_n^{\dbP} \geq t}{\sum}\int_t^{\tilde t} G^n(r, \dbP_{Y_r})dr + V(\tilde t, \dbP_{Y_{\tilde t-}}; \dbP_{\a^{t,l}_{\tilde t-}}) \Big\} ,\\
\dis V(\tilde t, m;l) \ge M[V] (\tilde t, m;l).
\ea\right.
\eea

First, recall (16) in \cite{LW}, 
\beaa
 \underset{i \in \cN}{\sum}\underset{j \in \cN}{\sum}\int_{\dbR^d}g_{ij}(t,x)p_{ij}m(dx,i)=\sum_{\t^{\dbP}_n \geq 0} G^n(t, m). 
\eeaa
For any $\dbP\in \cP(t, m;l)$, note that $m' := \dbP\circ (X_t, I_{\tilde t-})^{-1} \preceq m$, then
\beaa
 V(\tilde t, \dbP_{Y_{\tilde t-}};\dbP_{\a^{t,l}_{\tilde t-}})  - V(\tilde t, m;\dbP_{\a^{t,l}_{\tilde t-}})  &\le V(\tilde t, \dbP_{(X_{\tilde t}, I_{\tilde t-})};\dbP_{\a^{t,l}_{\tilde t-}})  - V(\tilde t, \dbP_{(X_t, I_{\tilde t-})};\dbP_{\a^{t,l}_{\tilde t-}}) + \underset{\t_n^{\dbP} \geq t}{\sum} G^n(\tilde t, m)\\
& \le C\dbE^\dbP[1+|X_t|]\sqrt{\tilde t-t} + \underset{\t_n^{\dbP} \geq t}{\sum} G^n(\tilde t, m)
\eeaa
Since $\dbP\in \cP(t, m;l)$ is arbitrary, by \reff{A-DPP} we have
\bea\label{dR2}
\begin{split}
V(t, m;l)  - V(\tilde t, m;l)  &\le \sup_{\dbP\in \cP(t, m;l)}  \Big\{\int_t^{\tilde t}  F(r, \dbP_{\a^{t,l}_r}, \dbP_{Y_r}) dr -\underset{\t_n^{\dbP} \geq t}{\sum}\int_t^{\tilde t} G^n(r, \dbP_{Y_r})dr + \int_t^{\tilde t} \dbQ V(\tilde t,m,\dbP_{\a^{t,l}_r})dr\\
& + C\dbE^\dbP[1+|X_t|]\sqrt{\tilde t-t} + \underset{\t_n^{\dbP} \geq t} {\sum}G^n(\tilde t, m) \Big\}  .
\end{split}
\eea

Next, for $m'\preceq m$, choose $\dbP\in \cP(t, m;l)$ s.t. $I_s = I_{t-}$, $t\le s< \tilde t$, then $\dbP\circ (X_t, I_{\tilde t-})^{-1} = m$. 
\beaa
 V(\tilde t, m,\dbP_{\a^{t,l}_{\tilde t-}}) - V(t, m;l) \le V(\tilde t, \dbP_{ (X_t, I_{\tilde t-})};\dbP_{\a^{t,l}_{\tilde t-}}) - V(\tilde t, \dbP_{ (X_{\tilde t}, I_{\tilde t-})};\dbP_{\a^{t,l}_{\tilde t-}}) - \int_t^{\tilde t}  F(r,\dbP_{\a^{t,l}_r}, \dbP_{Y_r}) dr .
 \eeaa
Thus, we have
\beaa
V(\tilde t, m;l)  - V(t, m;l)  \le  C\dbE^\dbP[1+|X_t|]\sqrt{\tilde t-t}- \int_t^{\tilde t}  F(r,\dbP_{\a^{t,l}_r},\dbP_{Y_r}) dr -\int_t^{\tilde t} \dbQ V(\tilde t,m,\dbP_{\a^{t,l}_r})dr.
\eeaa
The above together with \reff{dR2}, implies the desired regularity $\underset{\tilde t \downarrow t}{\lim} V(\tilde t,  m) = V(t,m)$ immediately.
\qed

\bibliographystyle{plain}

\end{document}